\documentclass[10pt]{article}
\usepackage{amssymb}
\usepackage{graphicx}
\usepackage{amsmath}
\usepackage{color}
\usepackage[a4paper]{geometry}

 \newtheorem{thm}{Theorem}[section]
 \newtheorem{cor}[thm]{Corollary}
 \newtheorem{lem}[thm]{Lemma}
 \newtheorem{prop}[thm]{Proposition}
 \newtheorem{defn}[thm]{Definition}
 \newtheorem{rem}[thm]{Remark}
\newtheorem{ex}[thm]{Example}
 \numberwithin{equation}{section}

\def\hh{\!\!\!\!}
\def\LS{Lebesgue-Stieltjes\ }
\def\PP{Proposition\ }
\def\RS{Riemann-Stieltjes\ }
\def\Proof{\noindent\emph{Proof.} }
\def\qed{\hfill $\Box$ \smallskip}

\begin{document}

\title{\textbf{Dependence of Solutions and Eigenvalues of Third Order Linear Measure Differential Equations on Measures}}
\author{Yixuan Liu$^1$, \qquad Guoliang Shi$^1$, \qquad Jun Yan$^{1}$\thanks{%
Corresponding author.}}
\date{}
\maketitle

\begin{center}
$^1$ School of Mathematics, Tianjin University, Tianjin, 300354, People's Republic of China

E-mail: \texttt{liuyixuan@tju.edu.cn, glshi@tju.edu.cn, jun.yan@tju.edu.cn}

\end{center}

\begin{abstract}
This paper deals with a complex third order linear measure differential equation
\begin{equation*}
i\mathrm{d}\left( y^{\prime }\right) ^{\bullet }+2iq\left(
x\right) y^{\prime }\mathrm{d}x+y\left( i\mathrm{d}q\left( x\right) +\mathrm{d}p\left(
x\right) \right) = \lambda y\mathrm{d}x
\end{equation*}
on a bounded interval with boundary conditions presenting a mixed aspect of the Dirichlet and the periodic problems. The dependence of eigenvalues on the {coefficients} $p$, $q$ is investigated. {We prove that the $n$-th eigenvalue is continuous in $p$, $q$ {when the norm topology of total variation and the weak$^*$ topology are considered}}. Moreover, the Fr\'{e}chet differentiability of the $n$-th eigenvalue {in} $p$, $q$ {with} the norm topology of total variation is also considered. To deduce these conclusions, we investigate the dependence of solutions of the above equation on the coefficients $p$, $q$ with different topologies and establish the counting lemma of eigenvalues according to the estimates of solutions.

\textbf{Mathematics Subject Classification (2010)}: 34A12; 34L40; 58C07

\textbf{Keywords}: Measure differential equation; Third order; Continuity; Eigenvalue
\end{abstract}

\section{Introduction}

In this paper, we consider the measure differential equation
\begin{equation}
i\mathrm{d}\left( y^{\prime }\right) ^{\bullet }+2iq\left(
x\right) y^{\prime }\mathrm{d}x+y\left( i\mathrm{d}q\left( x\right) +\mathrm{d}p\left(
x\right) \right) = \lambda y\mathrm{d}x, \qquad x\in {I:=[0,1]} \label{equation}
\end{equation}%
with the boundary
conditions%
\begin{equation*}
\left( BC\right)_{1} \ \ \left\{
\begin{array}{l}
y\left( 1\right) =0, \\
y^{\prime }\left( 1\right) =y^{\prime }\left( 0\right) , \\
\left( y^{\prime }\right) ^{\bullet }\left( 0\right) =0,%
\end{array}%
\right.  \label{bc+}
\end{equation*}
or the boundary conditions
\begin{equation*}
\left( BC\right)_{2} \ \ \left\{
\begin{array}{l}
y\left( 1\right) =0, \\
y^{\prime }\left( 1\right) =-y^{\prime }\left( 0\right) , \\
\left( y^{\prime }\right) ^{\bullet }\left( 0\right) =0,%
\end{array}%
\right.  \label{bc-}
\end{equation*}%
where $p$, $q\in {\mathcal{M}}_{0}(I,{\mathbb{R}})$ and $\lambda$ is a parameter in $\mathbb{C}$. Here{,} ${\mathcal{M}}_{0}(I,{\mathbb{R}})$ denotes the space of real-valued measures on ${I}$, which is the same as the dual space of the Banach space of continuous functions. The notation $y^{\prime }\left( x\right) $ stands for the classical derivative of $y\left( x\right) ${,} and $y^{\bullet }\left(
x\right) $ represents the generalized right-derivative of $y\left(
x\right) $ which will be defined precisely later (see Corollary \ref{solpp} (ii)).

Let $z\left( x\right) :=y^{\prime }\left( x\right) $, $w\left( x\right)
:=\left( y^{\prime }\right) ^{\bullet }\left( x\right) $, and then {the} equation (\ref%
{equation}) is equivalent to the following system %
\begin{equation}
\left\{
\begin{array}{l}
\mathrm{d}y(x)=z(x) \mathrm{d}x,\\
\mathrm{d}z(x)=w\left( x\right)\mathrm{d}x,\\
\mathrm{d}w(x)=-2q\left( x\right) z(x)\mathrm{d}x-y(x)\mathrm{d}\mu (x),
\end{array}%
\right. \label{ss}
\end{equation}
where $\mu (x)=q\left( x\right) -ip\left( x\right) +\lambda ix$. Therefore{,} using the facts of \LS integral, the solution of (\ref{equation}) with initial conditions
\begin{equation}
\left( y\left( 0\right) ,y^{\prime }\left( 0\right) ,\left( y^{\prime
}\right) ^{\bullet }\left( 0\right) \right) =\left( y\left( 0\right)
,z\left( 0\right) ,w\left( 0\right) \right) =(y_{0},z_{0},w_{0})\in {\mathbb{K}}^{3},\mathbb{K}=\mathbb{R} \mbox{ or }\mathbb{C} \label{initial}
\end{equation}%
is defined in Definition \ref{solss}. The imaginary {unit} $i$ in (\ref{equation}) indicates the solutions of this equation are usually complex-valued, even if $\lambda \in \mathbb{R}$; it is the reason why we call (\ref{equation}) a complex third order linear measure differential equation. It will be proved that the boundary value problem (\ref{equation})-$(BC)_{\xi}$, $\xi=1,2$, admits a real increasing sequence of eigenvalues
\begin{equation}
\Lambda _{\xi }\left(p,q\right) =\left\{ \lambda _{\xi,n}\left( p,q\right) , n\in\mathbb{Z}\right\} ,\xi=1,2,\label{L1}
\end{equation}
where $\mathbb{Z}=\left\{ 0,\pm 1,\pm 2,\cdots \right\} $ (see Lemma \ref{ev-real}) and the geometric multiplicity of each eigenvalue $\lambda _{\xi,n}$, $\xi=1,2$, is at most two (see Lemma \ref{doe}).

{Measure differential equations enable us to treat in a unified way both continuous and discrete systems, which} have attracted tremendous interest in the last decades. The researches on second and fourth order measure differential equations {can be} found in papers \cite{ali2003, Zh10, tes2013jdm, mg2013jde, mg2015pams, mg2016jde, die2018, dcds2018, jde2018} and the references therein. In contrast, third order measure differential equations have not yet been studied in the literature, and it is precisely the purpose of this paper to investigate the solutions and eigenvalues of the boundary value problems (\ref{equation})-(BC)$_{\xi}$, $\xi =1,2$.

{Note that in the special case $\left(\frac{\mathrm{d}p}{\mathrm{d}x},q\right) =:(u,v)\in {\mathcal{L}}^{2}(I,{\mathbb{R}})\times {\mathcal{H}}^{1}(I,{\mathbb{R}})$, the equation (\ref{equation}) reduces to the standard one
\begin{equation*}
L_{u,v}y= \lambda y,
\end{equation*}
where
\begin{equation*}
L_{u,v}:=iD^{3}+iDv+ivD+u,\qquad D:=\frac{\partial}{\partial x}.
\end{equation*}
}We emphasize that the operator $L_{u,v}$ {occurs in the inverse problem method of integration for the nonlinear evolution Boussinesq equation (see \cite{mckean1981cpam} for more considerations):
\begin{equation}
\frac{\partial ^{2} v}{\partial t^{2}}=\frac{\partial ^{2} }{\partial x^{2}}(\frac{4}{3}v^{2}+\frac{1}{3}\frac{\partial ^{2} v}{\partial x^{2}}),\quad   \frac{\partial u}{\partial x}= \frac{\partial v}{\partial t}.\label{bsq}
\end{equation}
Namely, (\ref{bsq}) is equivalent to the Lax equation $KL_{u,v}-L_{u,v}K=\dot{L}_{u,v}$, where $K=i(D^{2}+\frac{4}{3}v)$, and $\dot{L}_{u,v}$ denotes the derivative of ${L}_{u,v}$ with respect to $t$. Recently, the operator $L_{u,v}$} has attracted considerable attention {(see \cite{amour1999siam, amour2001siam, badanin2012jde, Badanin2013} and the references therein)}. In particular, for $\left(u,v\right) \in {\mathcal{L}}^{2}(I,{\mathbb{R}})\times {\mathcal{H}}^{1}(I,{\mathbb{R}})$, $v(0)=0$, Amour L {\cite{amour1999siam}} investigated the direct and inverse problems of operators $L_{u,v}$ on $I$ with the boundary conditions $\left( BC\right)_{\xi} $, $\xi=1,2$; the author discussed the multiplicities of eigenvalues, and then gave the {estimates} of solutions to deduce the counting lemma and estimates of the eigenvalues. In this paper, we first aim to generalize some results in \cite{amour1999siam} to the third order measure differential equation (\ref{equation}). More precisely, we show the estimates of solutions (see Theorem \ref{eopq}) of {the} equation (\ref{equation}), and then deduce the counting lemma (see Theorem \ref{countinglemma}) to illustrate the distribution, indexation and estimates (see Corollary \ref{estimate}) of eigenvalues, which is the first step towards the solution of the related inverse problem. {On the basis of} these results, we can characterize the dependence of the $n$-th eigenvalue $\lambda _{\xi,n}(p,q)$ on {the} {coefficients} $p$, $q$ as follows, which is the main result of this paper.

\begin{thm}
\label{scoe}
{Suppose $(p,q)\in {\mathcal{M}}_{0}(I,{\mathbb{R}})\times {\mathcal{%
M}}_{0}(I,{\mathbb{R}})$. }

{$(i)$ For any fixed $p \in {\mathcal{M}}_{0}(I,{\mathbb{R}})$, $\xi =1,2$, the eigenvalue $\lambda _{\xi ,n}\left( p,q\right) $ is continuous in $q\in {({\mathcal{M}}%
_{0}(I,{\mathbb{R}}),w^{\ast })}$.}

{$(ii)$ For any fixed $q \in {\mathcal{M}}_{0}(I,{\mathbb{R}})$, $\xi =1,2$, the eigenvalue $\lambda _{\xi ,n}\left( p,q\right) $ is continuous in $p\in ({\mathcal{M}}_{0}(I,{\mathbb{R}}),w^{\ast })$.}
\end{thm}

{Here, the symbol $({\mathcal{M}}_{0}(I,{\mathbb{R}}),w^{\ast })$ denotes the measure space with {the} weak$^{\ast}$ topology whose definition can be found in Section 2, and we use $({\mathcal{M}}_{0}(I,{\mathbb{R}}),\Vert \cdot \Vert _{\mathbf{V}})$ to denote the measure space with $\Vert \cdot \Vert _{\mathbf{V}}$-topology.} Note that {Theorem \ref{scoe}} indicates the eigenvalue $\lambda _{\xi ,n}\left( p,q\right) $ is also continuous in $p$, ${q}\in ({\mathcal{M}}_{0}(I,{\mathbb{R}}),\Vert \cdot \Vert _{\mathbf{V}})$ since the $w^{\ast }$-topology is weaker than the $\Vert \cdot \Vert _{\mathbf{V}}$-topology. In the {rest of this work}, we use
\begin{equation*}
E_{\xi ,n}\left( x,p,q\right) :=\frac{e(x,\lambda _{\xi ,n}\left(
p,q\right) ,p,q)}{\left( \int_{I}\left\vert e(x,\lambda _{\xi ,n}\left(
p,q\right) ,p,q)\right\vert ^{2}dx\right) ^{\frac{1}{2}}}
\end{equation*}
to denote the normalized eigenfunction corresponding to the simple eigenvalue $\lambda _{\xi ,n}\left(p,q\right)$, $\xi =1,2$. As a consequence of Theorem \ref{scoe}, the Fr\'{e}chet differentiability of eigenvalues with respect to $p$, $q \in ({\mathcal{M}}_{0}(I,{\mathbb{R}}),\Vert \cdot \Vert _{\mathbf{V}})$ can be obtained.

\begin{thm}
\label{frechet}
{$(i)$ Fix $q\in({\mathcal{M}}%
_{0}(I,{\mathbb{R}}),\Vert \cdot \Vert _{\mathbf{V}})$ and consider the eigenvalue $\lambda _{\xi ,n}\left( p,q\right) $, $\xi =1,2$ as a function of $p\in({\mathcal{M}}%
_{0}(I,{\mathbb{R}}),\Vert \cdot \Vert _{\mathbf{V}})$. Then for any $ p_{0}\in({\mathcal{M}}%
_{0}(I,{\mathbb{R}}),\Vert \cdot \Vert _{\mathbf{V}})$, there exists an integer $N_{1,p_{0}}>0$ such that $\lambda _{\xi ,n}\left( p,q\right) $, $\left\vert n\right\vert \geqslant N_{1,p_{0}}$, is continuously Fr\'{e}chet differentiable at $p_{0}$ and its Fr\'{e}chet derivative is given by
\begin{eqnarray*}
\partial _{p}\lambda _{\xi ,n}\left( p_{0},q\right)\!\!\!\!
&=&\!\!\!\!\left\vert E_{\xi ,n}\left( x,p_{0},q\right) \right\vert ^{2}\in ({%
\mathcal{C}}(I,{\mathbb{R}}),\Vert \cdot \Vert _{\infty })  \label{lp} \\
&\hookrightarrow &\!\!\!\!({\mathcal{C}}(I,{\mathbb{R}}),\Vert \cdot \Vert
_{\infty })^{\ast \ast }\cong ({\mathcal{M}}_{0}(I,{\mathbb{R}}),\Vert \cdot
\Vert _{\mathbf{V}})^{\ast },  \notag
\end{eqnarray*}
where ${\mathcal{C}}(I,{\mathbb{R}}):=\{g:I\rightarrow {\mathbb{R}}; g$ is
continuous on $I\}$, $\Vert g\Vert _{\infty }:=\underset{x\in I}{%
\sup }$ $\left\vert g\left( x\right) \right\vert $. }

{$(ii)$ Fix $p \in ({\mathcal{M}}%
_{0}(I,{\mathbb{R}}),\Vert \cdot \Vert _{\mathbf{V}})$ and consider the eigenvalue $\lambda _{\xi ,n}\left( p,q\right) $, $\xi =1,2$ as a function of $q\in({\mathcal{M}}%
_{0}(I,{\mathbb{R}}),\Vert \cdot \Vert _{\mathbf{V}})$. Then for any $q_{0}\in({\mathcal{M}}%
_{0}(I,{\mathbb{R}}),\Vert \cdot \Vert _{\mathbf{V}})$, there exists an integer $N_{2,q_{0}}>0$ such that $\lambda _{\xi ,n}\left( p,q\right) $, $\left\vert n\right\vert \geqslant N_{2,q_{0}}$, is continuously Fr\'{e}chet differentiable at $q_{0}$ and its Fr\'{e}chet derivative is given by
\begin{eqnarray*}
\partial _{q}\lambda _{\xi ,n}\left( p,q_{0}\right) \!\!\!\! &=&\!\!\!\! i
\left[ E_{\xi ,n}\left( x,p,q_{0}\right) ,\bar{E}_{\xi ,n}\left( x,p,q_{0}\right) %
\right] \in ({\mathcal{C}}(I,{\mathbb{R}}),\Vert \cdot \Vert _{\infty })
\label{lq} \\
&\hookrightarrow &\!\!\!\!({\mathcal{C}}(I,{\mathbb{R}}),\Vert \cdot \Vert
_{\infty })^{\ast \ast }\cong ({\mathcal{M}}_{0}(I,{\mathbb{R}}),\Vert \cdot
\Vert _{\mathbf{V}})^{\ast },  \notag
\end{eqnarray*}
where the operation $[\cdot,\cdot]$ is defined by $\left[ y,z\right] =zy^{\prime}-yz^{\prime}$.}
\end{thm}

{It is worth mentioning} that in \cite{mg2013jde}, Meng G and Zhang M considered the second order measure differential equation
\begin{equation}
\mathrm{d}y ^{\bullet }+\lambda y\mathrm{d}t+y\mathrm{d}\mu(t) = 0, \qquad t\in (0,1) \label{equation2}
\end{equation}%
with Neumann or Dirichlet boundary conditions, and investigated the dependence of eigenvalues on {the} measures $ \mu\in{\mathcal{M}}_{0}(I,{\mathbb{R}})$ with different topologies. Theorem \ref{scoe} and Theorem \ref{frechet} {generalize} the main results (Theorem 1.3 and Theorem 1.4) in \cite{mg2013jde}. Unfortunately, it seems to the authors that the approach in \cite{mg2013jde} cannot {apply} to this paper {directly because of the} major differences between the third order measure differential equation (\ref{equation}) and the second order measure differential equation (\ref{equation2}). For example, the solutions of {(\ref{equation})} are complex-valued and there exists the possibility of non-simple eigenvalues due to the coupled boundary conditions (BC)$_{\xi}$, $\xi=1,2$. Additionally, the eigenvalues of {the} boundary value problems (\ref{equation})-(BC)$_{\xi}$, $\xi=1,2$ are unbounded below and above. Nevertheless, we will propose a way to overcome these problems. In order to undertake the proofs, the dependence of solutions of {(\ref{equation})} on {the} measures $p$, $q$ with different topologies (see \PP \ref{wc}{, Remark \ref{uc1} and {\PP}\ref{uc2}}) and the counting lemma (see Theorem \ref{countinglemma}) are very crucial. It is also worth noting that the dependence of eigenvalues on the coefficients $p$, $q$ is of interest not only theoretically but also numerically. For classical Sturm-Liouville problems, Kong Q and Zettl A found that the numerical computation of the eigenvalue {is} based on the dependence of eigenvalues on the coefficients {(see \cite{kz96} and the references therein)}.

This paper is organized as follows. In Section 2, we introduce some basic definitions and useful properties of measures, \LS integral{,} and weak$^{\ast }$ topology{;} the existence and uniqueness of solutions {are also given}. Section 3 {investigates} the dependence of solutions on {the} measures $p$, $q$ with different topologies. Besides, we investigate the estimates of solutions and the analytic dependence of solutions on the spectral parameter $\lambda$. Finally, Section 4 provides the counting lemma to explain the distribution and asymptotic formulas of eigenvalues{;} the proof of the dependence of the $n$-th eigenvalue on {the} measures $p$, $q$ with different topologies is also given.

\section{Preliminaries}

\subsection{Measures, \LS Integral and Weak$^{\ast }$
Topology}

In this subsection, we briefly review some basic facts of
measures, different topologies of the measure space, \LS integral{,} and \RS integral. The detailed
theory can be founded in \cite{CvB00, DS58, Me98}.

Let ${\mathbb{N}}:=\left\{
1, 2,\cdots \right\} $ and ${\mathbb{K}}={\mathbb{R}}$ or ${\mathbb{C}}$. {Recall that $I=[0,1]$}. Then the space of (non-normalized) ${\mathbb{K}}$%
-value measures of $I$ is defined as
\begin{equation*}
{\mathcal{M}}(I,{\mathbb{K}}):=\{f:I\rightarrow {\mathbb{K}};f(0+)\
\exists ,f(x+)=f(x)\ \forall x\in (0,1),\mathbf{V}(f,I)<\infty \},
\end{equation*}%
where
\begin{equation*}
\mathbf{V}(f,I)\!:=\!\sup \!\left\{\!
\sum_{j=0}^{m-1}\!|f(x_{j+1})\!-\!f(x_{j})|\!:\!0\!=\!x_{0}\!<\!x_{1}\!<\!\cdots
\!<\!x_{m-1} \!<\!x_{m}\!=1\!,m\!\in\! {\mathbb{N}}\!\right\}
\end{equation*}%
is the \textit{total variation}\ of $f$ over $I$ and {for any $x\in \lbrack 0,1)$, $f(x+):=\underset{t\rightarrow
x+}\lim f(t)$ denotes the right-limit}. Note that ${%
\mathcal{M}}(I,{\mathbb{K}})$ is a Banach space with the norm $\Vert
f\Vert _{\mathbf{V}}=\mathbf{V}(f,I)+\left\vert f(0)\right\vert $. The total variation of $f$ over any subinterval $I_{0}$ (closed, open or semi open) is also well-defined. For example, if $%
I_{0}=(a,b]\subset I$, the total variation is
\begin{equation*}
\mathbf{V}(f,I_{0})\!:=\!\sup\! \left\{\!
\sum_{j=0}^{m-1}\!|f(x_{j+1})\!-\!\!f(x_{j})|:a<x_{0}<x_{1}<\cdots
<x_{m-1} <x_{m}=b, m\in {\mathbb{N}}\right\} .
\end{equation*}%
{F}or any $x\in (0,1)$, $f(x+)=f(x)$, {thus} we obtain that for each $x_{0}\in (0,1)$,%
\begin{equation*}
\lim_{x\rightarrow x_{0}+}\mathbf{V}(f,[x_{0},x])=\lim_{x\rightarrow x_{0}+}%
\mathbf{V}(f,(x_{0},x])=0.
\end{equation*}%

The space of (normalized) ${\mathbb{K}}$-valued measures is
\begin{equation*}
{\mathcal{M}}_{0}(I,{\mathbb{K}}):=\left\{ f\in {\mathcal{M}}(I,{\mathbb{K}}%
):f(0)=0\right\}{,}
\end{equation*}%
and the normalization condition for $f\in {\mathcal{M}}_{0}(I,{\mathbb{K}})$ is $%
f(0)=0$. Hence{,} {$f(0+)\neq 0$ is possible and $\mathbf{V}%
(f,I)=\Vert f\Vert _{\mathbf{V}}$}. The topology induced by the norm $\Vert \cdot \Vert _{\mathbf{V}}$ is called the strong topology ({$\Vert \cdot \Vert _{\mathbf{V}}$-topology}) of ${\mathcal{M}}_{0}(I,{\mathbb{K}})$. According to the Riesz representation theorem, $({\mathcal{M}}_{0}(I,{\mathbb{K}}%
),\Vert \cdot \Vert _{\mathbf{V}})$ is {identical to} the dual space of the
Banach space $({\mathcal{C}}(I,{\mathbb{K}}),\Vert \cdot \Vert _{\infty })$,
where ${\mathcal{C}}(I,{\mathbb{K}}):=\{g:I\rightarrow {\mathbb{K}}; g$ is
continuous on $I\}$, $\Vert g\Vert _{\infty }:=\underset{x\in I}{%
\sup }$ $\left\vert g\left( x\right) \right\vert $.

In fact, {any} $f\in ({%
\mathcal{M}}_{0}(I,{\mathbb{K}}),\Vert \cdot \Vert _{\mathbf{V}})$ defines $%
f^{\ast }\in ({\mathcal{C}}(I,{\mathbb{K}}),\Vert \cdot \Vert _{\infty
})^{\ast }$ by
\begin{equation}
f^{\ast }(g)=\int_{I}g(t)\mathrm{d}f(t),\qquad g\in {\mathcal{C}}(I,{%
\mathbb{K}}),  \label{dual}
\end{equation}%
which refers to the \RS integral. Moreover, one has
\begin{equation*}
\Vert f\Vert _{\mathbf{V}}=\mathbf{V}(f,I)=\sup \left\{
\mbox{$\int_{I}
g\mathrm{d}f: g\in {\mathcal C}(I,{\mathbb K}),  \Vert g\Vert _{\infty}=1$}\right\} .
\end{equation*}%
From the duality relation (\ref{dual}), we define the following weak$^{\ast }$ topology $w^{\ast }$ of ${\mathcal{M}}_{0}(I,{\mathbb{K}})$.

\begin{defn}
\label{ws} {For $f_{0}, f_{m}\in {\mathcal{M}}_{0}(I,{\mathbb{K}})$,
$m\in {\mathbb{N}}$, we say $f_{m}$ is weakly$^{\ast }$ convergent to $f_{0}$ {as $m\rightarrow \infty$}{,} if and only if for each $g\in {\mathcal{C}}(I,{\mathbb{K}})$,
\begin{equation*}
\lim_{m\rightarrow \infty }\int_{I}g\mathrm{d}f_{m}=\int_{I}g\mathrm{d}f_{0}\mathrm{.}
\end{equation*}}
\end{defn}
Apparently, {the following example illustrates} the weak$^{\ast }$ topology is weaker than $\Vert \cdot \Vert_{\mathbf{V}}${-topology}.

\begin{ex}
\label{ex1}{For $a\in \left( 0,1 \right]$, let
\begin{equation*}
\delta _{a}(x):=\left\{
\begin{array}{l}
0\mbox{ for }x\in \lbrack 0,a), \\
1\mbox{ for }x\in \lbrack a,1],%
\end{array}%
\right.
\end{equation*}
and
\begin{equation*}
\delta _{0}(x):=\left\{
\begin{array}{l}
0\mbox{ for }x=0, \\
1\mbox{ for }x\in (0,1].%
\end{array}%
\right.
\end{equation*}
{For any $g \in \mathcal{C}(I,\mathbb{R})$,} we have $\int_{I}g\mathrm{d}\delta _{a}=g(a)\rightarrow g(0)=\int_{I}g\mathrm{d}\delta _{0}$ as $a\rightarrow 0$, {{i.e.}, $\delta _{a} \rightarrow \delta _{0}$ in $(\mathcal{M}_{0}(I, \mathbb{K}),w^{\ast})$ as $a\rightarrow 0$. However,} for any $a \in \left(0,1\right]$, $ \Vert\delta _{a}-\delta _{0} \Vert_{\mathbf{V}}=2$.}
\end{ex}

In \cite{HT09,Tv02}, another topology induced by the supremum norm $\Vert \cdot \Vert _{\infty }$ is also used for ${%
\mathcal{M}}_{0}(I,{\mathbb{K}})$. As $\Vert f\Vert _{\infty }\leqslant
\Vert f\Vert _{\mathbf{V}}$ for all $f\in {\mathcal{M}}_{0}(I,{\mathbb{K}})$%
, one sees that $\Vert \cdot \Vert _{\infty }$ is also weaker than $\Vert
\cdot \Vert _{\mathbf{V}}$. Moreover, we obtain the following relations for the weak$^{\ast }$ topology and the topology induced by the norm $\Vert \cdot \Vert _{\infty }$.

\begin{lem}
\label{wt-rel} {One has
\begin{equation}  \label{TO1}
f_m \to f_0 \mbox{ in } ({\mathcal{M}}_0(I,{\mathbb{K}}),\|\cdot\|_\infty)%
\not \Longrightarrow f_m \to f_0 \mbox{ in } ({\mathcal{M}}_0(I,{\mathbb{K}}),w^{\ast }),
\end{equation}
\begin{equation}  \label{TO2}
f_m \to f_0 \mbox{ in } ({\mathcal{M}}_0(I,{\mathbb{K}}),w^{\ast }) \not \Longrightarrow f_m
\to f_0\mbox{ in } ({\mathcal{M}}_0(I,{\mathbb{K}}),\|\cdot\|_\infty).
\end{equation}
}
\end{lem}

\Proof Let $f_m(x)=\frac {1}{m}\sin (2\pi m^{2}x) \in {\mathcal{M}}_0(I,{\mathbb{R}})$, $m \in \mathbb{N}$, and $f_0(x)=0$, then $f_{m} \rightarrow f_{0} $ in $({\mathcal{M}}_0(I,{\mathbb{K}}),\|\cdot\|_\infty)$. Since $\int_{I}\mathrm{d}f _{m}=4m \rightarrow \infty $, we know the relation (\ref{TO1}) holds. From Example \ref{ex1}, one has $\|\delta _{a}-\delta _{0}\|_\infty =1 \neq 0$ holds for $a \in \left(0,1\right]$, and thus we obtain the relation (\ref{TO2}). \qed

Given $f\in {\mathcal{M}}_{0}(I,{\mathbb{K}})$ and $g\in {\mathcal{C}}(I,{%
\mathbb{K}})$, for any subinterval $I_{0}\subset I$, the \LS
integral $\int_{I_{0}}g \mathrm{d}f$ is also defined. Due to the possible
jump of a measure $f(x)$ at $x=0$, one has
\begin{equation*}
\int_{\lbrack 0,b]}g\mathrm{d}f=g(0)f(0+)+\int_{(0,b]}g \mathrm{d}f,\qquad
b\in (0,1]{,}  \label{jf1}
\end{equation*}%
{{{i.e.}}}, $\int_{[0,b]}g \mathrm{d}f$ and $\int_{(0,b]}g\mathrm{d}f$ may
differ. If $I_{0}$ has the form $(a,b)$, $(a,b]$, where $0\leqslant
a<b\leqslant 1$, or the form $[0,b)$, $[0,b]$, where $0<b\leqslant 1$, one
has the following basic inequality
\begin{equation*}
\left\vert \int_{I_{0}}g \mathrm{d}f\right\vert \leqslant \Vert g\Vert
_{\infty ,I_{0}}\cdot \mathbf{V}(f,I_{0}),\qquad \mbox{where }\Vert g\Vert
_{\infty ,I_{0}}:=\sup_{t\in I_{0}}|g(t)|.
\end{equation*}
For real measures, we have the following lemmas.
\begin{lem}
\label{Vmu1}{For $f\in {\mathcal{M}}_{0}(I,{\mathbb{R}})$, let
\begin{equation*}
\mathbf{V}_{f}\left( x\right) :=\left\{
\begin{array}{ll}
0, & x=0,\\
\mathbf{V}(f,\left( 0,x\right] ),& x \in (0,1],
\end{array}
\right.
\end{equation*}
then we have $\mathbf{V}_{f}\left( x\right) \in {\mathcal{M}}_{0}(I,{\mathbb{R}})$ and
\begin{equation}
\left\vert \int_{\lbrack a,b]}g(x) \mathrm{d}f(x)\right\vert \leqslant
\int_{\lbrack a,b]}\vert g(x)\vert \mathrm{d}\mathbf{V}_{f}\left( x\right) \qquad
\forall g\in {\mathcal{C}}(I,{\mathbb{K}}), [a,b]\subset I.  \label{i3}
\end{equation}}
\end{lem}

\Proof See \cite[p. 321]{xdx}. \qed

\begin{lem}
\label{ubp}
{Suppose the sequence $\{f_{m}\}$ converges to $f_{0}$ in $({\mathcal{M}}_0(I,{\mathbb{K}}),w^*)$, then there exists a constant $C^{\ast}_{f_0}>0$ such that $\underset{m\in {\mathbb{N}_{0}}}\sup\Vert f_{m}\Vert _{\mathbf{V}}\leqslant C_{f_{0}}^{\ast}$. }
\end{lem}
\Proof {Due to the fact that weak$^{\ast }$ {convergence} implies boundedness}, this lemma can be proved. \qed

\subsection{Notation, Existence and Uniqueness of Solutions}

In the following, we give some basic facts on the solutions of (\ref{equation}), where $p$, $q\in {\mathcal{M}}_{0}(I,{\mathbb{K}})$, $\lambda \in {\mathbb{C}%
}$. Due to the equivalence between {the} equation (\ref{equation}) and {the} system (\ref{ss}), the solution of (\ref{equation}) with initial conditions (\ref{initial}) is defined as follows.

\begin{defn}
\label{solss} {For $p$, $q \in {\mathcal{M}}_{0}(I,{\mathbb{K}})$, $\lambda \in \mathbb{C}$, $(y_{0},z_{0},w_{0})\in {\mathbb{K}}^{3}$, {a} function $y(x)$ is {a} solution of {the} initial value problem $(\ref{equation})$, $(\ref{initial})$ if it satisfies that }

$(i)$ {$y\in C^{1}(I,{\mathbb{C}}):=\left\{ f:I\rightarrow {%
\mathbb{C}}; f\text{ is continuously differentiable on }I\right\}
$, and }

$(ii)$ {there exist functions $z, w:I\rightarrow {\mathbb{C}}$ such that
\begin{eqnarray*}
&&y(x)=y_{0}+\int_{[0,x]}z(t)\mathrm{d}t,\qquad x\in I\text{,}  \label{ryz} \\
&&z(x)=z_{0}+\int_{[0,x]}w(t)\mathrm{d}t,\qquad x\in I,
\label{rzw} \\
&&w(x)=\left\{
\begin{array}{ll}
w_{0}, & x=0, \\
w_{0}-\int_{[0,x]}2q(t)z(t)\mathrm{d}t-\int_{[0,x]}y(t) \mathrm{d}\mu (t), & x\in
(0,1].%
\end{array}%
\right.  \label{rwy}
\end{eqnarray*}%
}
\end{defn}

{{The solution $y$ is defined via fixed point equations, and we can prove the existence and
uniqueness of the solution by many methods, one of which is based on the
Kurzweil-Stieltjes integral{, see \cite{HT09}}.

\begin{prop}
\label{e-u}{For each $(y_{0},z_{0},w_{0})\in {\mathbb{K}}^{3}$, the initial value problem $(\ref{equation})$, $(\ref{initial})$ has the unique solution $y(x)$ on $I$.}
\end{prop}

Since the solution $y$ is continuous differentiable on $I$, {one has} $ z\in {\mathcal{C}}(I,{\mathbb{C}})$, $w\in {\mathcal{M}}(I,{\mathbb{C}})\subset {\mathcal{L}}^{1}(I,{\mathbb{C}})$. {If we use $y^{\prime}$, $\left(y^{\prime}\right)^{\bullet}$ to denote $z$, $w$, respectively, then we have
\begin{eqnarray}
&&y(x)=y_{0}+\int_{[0,x]}y^{\prime }(t)\mathrm{d}t,\qquad x\in I\text{,}  \label{ryz} \\
&&y^{\prime }(x)=z_{0}+\int_{[0,x]}\left( y^{\prime
}\right) ^{\bullet }(t)\mathrm{d}t,\qquad x\in I,
\label{rzw} \\
&&\left( y^{\prime
}\right) ^{\bullet }(x)=\left\{
\begin{array}{ll}
w_{0}, & x=0, \\
w_{0}-\int_{[0,x]}2q(t)y^{\prime }(t)\mathrm{d}t-\int_{[0,x]}y(t) \mathrm{d}\mu (t), & x\in
(0,1].%
\end{array}%
\right.  \label{rwy}
\end{eqnarray}
}}}According to the property of
Lebesgue integral and \LS integral, we obtain the following
corollary.

\begin{cor}
\label{solpp}
{$(i)$ {There holds}}%
\begin{equation*}
\int_{\lbrack x_{1},x_{2}]}\left( y^{\prime }\right) ^{\bullet }(t) \mathrm{%
d}t=\int_{(x_{1},x_{2}]}\left( y^{\prime }\right) ^{\bullet }(t) \mathrm{d}%
t=y^{\prime }(x_{2})-y^{\prime }(x_{1}),\qquad 0\leqslant x_{1}\leqslant
x_{2}\leqslant 1.  \label{pro1}
\end{equation*}

{$(ii)$ $y^{\prime }$ is the classical derivative of $y$ with respect to $x$ on $I$, and $\left( y^{\prime }\right)^{\bullet }(x_{0})$ is the classical right-derivative at any point $x_{0}\in (0,1)$, i.e.,
\begin{equation*}
\left( y^{\prime }\right) ^{\bullet }(x_{0})=\lim_{x\rightarrow x_{0}+}\frac{%
y^{\prime }(x)-y^{\prime }(x_{0})}{x-x_{0}}.  \label{der10}
\end{equation*}}

{$(iii)$ Actually, $y^{\prime }$ is absolutely continuous on $I${. Hence,} the following identity
\begin{equation*}
\left( y^{\prime }\right) ^{\bullet }(x_{0})=y^{\prime \prime
}(x_{0}):=\lim_{x\rightarrow x_{0}}\frac{y^{\prime }(x)-y^{\prime }(x_{0})}{%
x-x_{0}} \label{dot1}
\end{equation*}
holds for Lebesgue-a.e. $x_{0}\in I$.}
\end{cor}

\Proof The proof is similar to {that of} \cite[Corollary 3.4]%
{mg2013jde}. \qed

In this paper, we use $%
y_{1}(x,\lambda ,p,q)$, $y_{2}(x,\lambda ,p,q)$, $y_{3}(x,\lambda ,p,q)$ to
denote the solutions of (\ref{equation}) satisfying the initial conditions%
\begin{equation*}
\left(
\begin{array}{ccc}
y_{1}(0,\lambda ,p,q) & y_{2}(0,\lambda ,p,q) & y_{3}(0,\lambda ,p,q) \\
y_{1}^{\prime }(0,\lambda ,p,q) & y_{2}^{\prime }(0,\lambda ,p,q) &
y_{3}^{\prime }(0,\lambda ,p,q) \\
\left( y_{1}^{\prime }\right) ^{\bullet }(0,\lambda ,p,q) & \left(
y_{2}^{\prime }\right) ^{\bullet }(0,\lambda ,p,q) & \left( y_{3}^{\prime
}\right) ^{\bullet }(0,\lambda ,p,q)%
\end{array}%
\right) =I_{3}:=\left(
\begin{array}{ccc}
1 & 0 & 0 \\
0 & 1 & 0 \\
0 & 0 & 1%
\end{array}%
\right).
\end{equation*}%
 {Denote}%
\begin{equation*}
N_{p,q}(x):=\left(
\begin{array}{ccc}
y_{1}(x,\lambda ,p,q) & y_{2}(x,\lambda ,p,q) & y_{3}(x,\lambda ,p,q) \\
y_{1}^{\prime }(x,\lambda ,p,q) & y_{2}^{\prime }(x,\lambda ,p,q) &
y_{3}^{\prime }(x,\lambda ,p,q) \\
\left( y_{1}^{\prime }\right) ^{\bullet }(x,\lambda ,p,q) & \left(
y_{2}^{\prime }\right) ^{\bullet }(x,\lambda ,p,q) & \left( y_{3}^{\prime
}\right) ^{\bullet }(x,\lambda ,p,q)%
\end{array}%
\right) , x\in I.  \label{npq}
\end{equation*}%
Then due to \PP \ref{e-u}, the solution
of {the} initial value problem (\ref{equation}), (\ref{initial}) can be denoted by
\begin{equation*}
\left(
\begin{array}{c}
y(x,\lambda ,p,q) \\
y^{\prime }(x,\lambda ,p,q) \\
\left( y^{\prime }\right) ^{\bullet }(x,\lambda ,p,q)%
\end{array}%
\right) =N_{p,q}(x)\left(
\begin{array}{c}
y_{0} \\
z_{0} \\
w_{0}%
\end{array}%
\right) , \left(
\begin{array}{c}
y_{0} \\
z_{0} \\
w_{0}%
\end{array}%
\right) \in {\mathbb{K}}^{3}.
\end{equation*}%

\begin{rem}
\label{ll}
{Since $N_{p,q}(0)=I_{3}$, the equality
\begin{equation*}
\det N_{p,q}(x)\equiv 1,\qquad x\in I  \label{ll1}
\end{equation*}
can be deduced by the same methods as those in $\cite{mg2013jde, Mi83}$.}
\end{rem}

\begin{lem}
\label{voc0} {The unique solution $(y(x),y^{\prime }(x),\left(
y^{\prime }\right) ^{\bullet }(x))$ of the third order inhomogeneous differential equation
\begin{eqnarray}
 i\mathrm{d}\left( y^{\prime }\right) ^{\bullet }+2iq\left( x\right)
y^{\prime }\mathrm{d}x+y\left( i\mathrm{d}q\left( x\right) +\mathrm{d}p\left( x\right)
-\lambda \mathrm{d}x\right) =ih(x) \mathrm{d}\nu (x), && \label{ided} \\
 p,q,\nu \in {\mathcal{M}}_{0}(I,{\mathbb{K}}),h\in {\mathcal{C}}(I,{\mathbb{K}})&&   \notag
\end{eqnarray}%
satisfying the initial conditions $(\ref{initial})$ is given by the {variation of constants} formula %
\begin{equation*}
\left(
\begin{array}{c}
y(x) \\
y^{\prime }(x) \\
\left( y^{\prime }\right) ^{\bullet }(x)%
\end{array}%
\right) =N_{p,q}(x)\left( \left(
\begin{array}{c}
y_{0} \\
z_{0} \\
w_{0}%
\end{array}%
\right) +\int_{[0,x]}N_{p,q}^{-1}(t)\left(
\begin{array}{c}
0 \\
0 \\
h(t)%
\end{array}%
\right)  \mathrm{d}\nu (t)\right), x\in (0,1].  \label{voc}
\end{equation*}%
Here{,} $N_{p,q}^{-1}(t)$ is the inverse of $N_{p,q}(t)$.}
\end{lem}
\Proof See \cite{Mi83}. \qed

\section{The Properties of Solutions of Measure Differential
Equation}
In this section, we investigate the dependence of the solution $y(x,\lambda ,p,q)$ and its derivatives $y^{\prime }(x,\lambda ,p,q)$, $\left( y^{\prime }\right) ^{\bullet}(x,\lambda ,p,q)$ on {the} measures $p$, $q\in {\mathcal{M}}_{0}(I,{\mathbb{K}})$ with different topologies. And then we give estimates of solutions and the analytic dependence of solutions on the spectral parameter $\lambda$ when $p$, $q\in {\mathcal{M}}_{0}(I,{\mathbb{R}})$.
\subsection{Dependence of Solutions on Measures $p$, $q$}

Firstly, we discuss the dependence of $y(x,\lambda,p,q) $, $y^{\prime }(x,\lambda ,p,q)$, $\left( y^{\prime }\right) ^{\bullet}(x,\lambda ,p,q)$ on {the} measure{s} $p$, ${q}\in ({\mathcal{M}}_{0}(I,{\mathbb{K}}),w^{\ast })$, which will be used in the proof of Theorem \ref{scoe}. The norm of $y\in C^{1}(I,{\mathbb{C}})$ is defined by $\Vert y\Vert _{{\mathcal{C}}^{1}}:=\Vert y\Vert _{\infty }+\Vert y^{\prime }\Vert _{\infty }$.

\begin{prop}
\label{wc}{{$(i)$} For any $\lambda \in {\mathbb{C}}$, the following mappings for the solution of {the} initial value problem $(\ref{equation})$, $(\ref{initial})$ are continuous,
\begin{eqnarray}
&&({\mathcal{M}}_{0}(I,{\mathbb{K}}),w^{\ast })\rightarrow \left( C^{1}(I,{%
\mathbb{C}}),\Vert \cdot \Vert _{{\mathcal{C}}^{1}}\right) ,\qquad
p\rightarrow y(\cdot ,\lambda ,p,q),  \label{py} \\
&&({\mathcal{M}}_{0}(I,{\mathbb{K}}),w^{\ast })\rightarrow \left( {\mathcal{M%
}}(I,{\mathbb{C}}),w^{\ast }\right) ,\qquad \ \ \ \ p\rightarrow \left(
y^{\prime }\right) ^{\bullet }(\cdot ,\lambda ,p,q).  \label{coy2}
\end{eqnarray}%
In particular, the following functional is continuous,
\begin{equation}
({\mathcal{M}}_{0}(I,{\mathbb{K}}),w^{\ast })\rightarrow {\mathbb{C}},\qquad
p\rightarrow \left( y^{\prime }\right) ^{\bullet }(1,\lambda ,p,q).
\label{coy21}
\end{equation}
{$(ii)$ For any $\lambda \in {\mathbb{C}}$,  the following mappings for the solution of {the} initial value problem $(\ref{equation})$, $(\ref{initial})$ are continuous,
\begin{eqnarray}
&&({\mathcal{M}}_{0}(I,{\mathbb{K}}),w^{\ast })\rightarrow \left( C^{1}(I,{\mathbb{C}}),\Vert \cdot \Vert _{{\mathcal{C}}^{1}}\right) ,\qquad q\rightarrow y(\cdot ,\lambda ,p,q), \label{yq}\\
&&({\mathcal{M}}_{0}(I,{\mathbb{K}}),w^{\ast })\rightarrow \left( {\mathcal{M%
}}(I,{\mathbb{C}}),w^{\ast }\right) ,\qquad\ \ \ \  q\rightarrow \left(
y^{\prime }\right) ^{\bullet }(\cdot ,\lambda ,p,q), \label{y2q}\\
&&({\mathcal{M}}_{0}(I,{\mathbb{K}}),w^{\ast })\rightarrow {\mathbb{C}}%
,\quad\qquad\qquad\qquad\quad\ \  q\rightarrow \left( y^{\prime }\right)^{\bullet }(1,\lambda ,p,q). \label{y1q}
\end{eqnarray}}}
\end{prop}

Before proving this proposition, we introduce some notations and a useful lemma as follows.

Assume that the sequence $\{p_{m}\}_{m\in \mathbb{N}}$ converges to $p_{0}$ in $({%
\mathcal{M}}_{0}(I,{\mathbb{K}}),w^{\ast })$. Let
\begin{eqnarray*}
&&y_{m}(x):=y(x,\lambda ,p_{m},q), \\
&&z_{m}(x):=y_{m}^{\prime }(x)=y^{\prime }(x,\lambda ,p_{m},q), \\
&&w_{m}(x):=\left( y_{m}^{\prime }\right) ^{\bullet }(x)=\left( y^{\prime
}\right) ^{\bullet }(x,\lambda ,p_{m},q),\qquad m\in {\mathbb{N}_{0}}:=\mathbb{N}\cup \{0\}.
\end{eqnarray*}%
We define functions $F, G:[0,1]^{2}\rightarrow {\mathbb{R}}$ by%
\begin{eqnarray*}
&&F(x,t):=\left\{
\begin{array}{ll}
2\left( x-t\right) & \mbox{ for }0\leqslant t\leqslant x\leqslant 1, \\
0 & \mbox{ for }0\leqslant x<t\leqslant 1,%
\end{array}%
\right. \\
&&G(x,t):=\left\{
\begin{array}{ll}
\frac{1}{2}\left( x-t\right) ^{2} & \mbox{ for }0\leqslant t\leqslant
x\leqslant 1, \\
0 & \mbox{ for }0\leqslant x<t\leqslant 1,%
\end{array}%
\right.
\end{eqnarray*}%
then we obtain $F$, $G$, $\frac{\partial }{\partial x}G\in {\mathcal{C}}(I^{2},{%
\mathbb{R}})$. For any $\lambda \in \mathbb{C}$, $p\in ({\mathcal{M}}_{0}(I,{\mathbb{K}}),w^{\ast })$, $y\in \left( C^{1}(I,{\mathbb{C}}),\Vert \cdot \Vert _{{\mathcal{C}}%
^{1}}\right) $,
{{{using the integration by parts formula for (\ref{rwy}) and the fact $q(0)=0$, we obtain
\begin{eqnarray*}
\left( y^{\prime}\right) ^{\bullet }(x)\!\!\!\!&=&\!\!\!\!w_{0}-2q(t)y(t)\vert_{t=0}^{t=x+}+\int_{[0,x]}y(t) \mathrm{d}(2q(t))-\int_{[0,x]}y(t) \mathrm{d}\mu (t)\\
\!\!\!\!&=&\!\!\!\!w_{0}-2q(x)y(x)+\int_{[0,x]}y(t)\mathrm{d}\tilde{\mu} (t),\quad x\in(0,1],
\end{eqnarray*}
where $\tilde{\mu}(t)=q(t)+ip(t)-\lambda it$. Therefore, we have
\begin{equation}
\left( y^{\prime }\right) ^{\bullet }(x)=\left\{
\begin{array}{ll}
w_{0}, & x=0, \\
w_{0}-2q(x)y(x)+\int_{[0,x]}y(t) \mathrm{d}\tilde{\mu}(t), & x\in (0,1].%
\end{array}%
\right.\label{ypb}
\end{equation}
Substitution of (\ref{ypb}) into (\ref{rzw}) yields
\begin{equation*}
y^{\prime }(x)=z_{0}+w_{0}x-\int_{[0,x]}2q(t)y(t)\mathrm{d}t+\int_{[0,x]}\int_{[0,t]}y(s) \mathrm{d}\tilde{\mu} (s)\mathrm{d}t, \quad x\in I.
\end{equation*}}
{Exchanging the order of integration in the double integral}, we find
\begin{equation}
y^{\prime}(x)=z_{0}+w_{0}x-\int_{[0,x]}2q(t)y(t)\mathrm{d}t+\int_{[0,x]}(x-t)y(t)\mathrm{d}\tilde{\mu} (t), \quad x\in I.\label{yprime}
\end{equation}}
}
{Substituting (\ref{yprime}) into (\ref{ryz}) and {exchanging the order of integration in the double integral} yield
\begin{eqnarray*}
y(x)\!\!\!\!&=&\!\!\!\!y_{0}+{\int_{[0,x]}[z_{0}+w_{0}t-\int_{[0,t]}2q(s)y(s)\mathrm{d}s+\int_{[0,t]}(t-s)y(s)\mathrm{d}\tilde{\mu} (s)]\mathrm{d}t}\notag\\
\!\!\!\!&=&\!\!\!\!\tilde{y}_{0}(x)-\int_{I}F(x,t)q\left( t\right) y(t) \mathrm{d}t+\int_{I}G(x,t)y(t) \mathrm{d}\tilde{\mu}(t), \quad x\in I,
\end{eqnarray*}
where $\tilde{y}_{0}(x)=y_{0}+z_{0}x+\frac{1}{2}w_{0}x^{2}$. }{Denote}
\begin{equation}
{\mathcal{Z}}\left( p,y\right) (x):=\int_{I}F(x,t)q\left( t\right) y(t) %
\mathrm{d}t-\int_{I}G(x,t)y(t) \mathrm{d}\tilde{\mu}(t),\ x\in I. \label{zpy}
\end{equation}%
Then a function $y\in {\mathcal{C}}^{1}(I,{\mathbb{C}})$ is a solution of {the} initial value problem (\ref{equation}), (\ref{initial}){,} if and only if it {satisfies}
\begin{equation}
y(x)=\tilde{y}_{0}(x)-{\mathcal{Z}}\left( p,y\right) (x).  \label{fix90}
\end{equation}
\begin{lem}
\label{rc}
{For any $\lambda \in {\mathbb{C}}$, $(y_{0},z_{0},w_{0})\in \mathbb{K}^{3}$,
the sequence $\{y_{m}\}_{m\in {\mathbb{N}_{0}}}$ is relatively
compact in $({\mathcal{C}}^{1}(I,{\mathbb{C}}),\Vert \cdot \Vert _{{%
\mathcal{C}}^{1}})$.}
\end{lem}
\Proof The proof of this lemma consists of three steps.

\textit{Step 1.} We need to verify that the sequence $\{y_{m}\}_{m\in {\mathbb{N}_{0}}}$ is uniformly
bounded.

Since the sequence $\{p_{m}\}_{m\in \mathbb{N}}$ converges to $p_{0}$ in $({%
\mathcal{M}}_{0}(I,{\mathbb{K}}),w^{\ast })$, it follows from Lemma \ref{ubp} that
\begin{equation*}
\underset{m\in {\mathbb{N}_{0}}}\sup\Vert p_{m}\Vert _{\infty }\leqslant \sup_{m\in {%
\mathbb{N}_{0}}}\Vert p_{m}\Vert _{\mathbf{V}}<\infty .
\end{equation*}%
{According to the integral equations (\ref{ryz}), (\ref{rzw}) and the definitions of $y_{m}(x)$, $z_{m}(x)$ and $w_{m}(x)$, one has
\begin{eqnarray}
|y_{m}(x)|\!\!\!\! &\leqslant &\!\!\!\!|y_{0}|+\int_{[0,x]}|z_{m}(t)| %
\mathrm{d}t,\qquad x\in I,  \label{ymzm} \\
|z_{m}(x)|\!\!\!\! &\leqslant &\!\!\!\!|z_{0}|+\int_{[0,x]}|w_{m}(t)| %
\mathrm{d}t,\qquad x\in I.  \label{zmwm}
\end{eqnarray}
From (\ref{ypb}), we find
\begin{equation*}
w_{m}(x)=\left\{
\begin{array}{ll}
w_{0}, & x=0, \\
w_{0}-2q(x)y_{m}(x)+\int_{[0,x]}y_{m}(t) \mathrm{d}(q(t)+ip_{m}(t)-\lambda it), & x\in
(0,1].%
\end{array}%
\right.
\end{equation*}
Hence,
\begin{eqnarray}
\vert w_{m}(x)\vert\!\!\!\! &\leqslant &\!\!\!\!\vert w_{0}\vert+2\Vert q\Vert _{\infty }\underset{t\in \lbrack 0,x]}\max \left\vert y_{m}(t)\right\vert+\underset{t\in \lbrack 0,x]}\max \left\vert y_{m}(t)\right\vert\left\vert \int_{[0,x]}\mathrm{d}(q(t)+ip_{m}(t)-\lambda it)\right\vert \notag\\
\!\!\!\! &\leqslant &\!\!\!\!\vert w_{0}\vert+2\Vert q\Vert _{\infty}\hat{y}_{m}(x)+\hat{y}_{m}(x)(\Vert q\Vert _{\mathbf{V}}+\Vert p_{m}\Vert _{\mathbf{V}}+\left\vert \lambda \right\vert) \notag\\
&\leqslant &\!\!\!\!|w_{0}|+C_{1}\hat{y}_{m}(x),\quad x\in I, \label{wmym}
\end{eqnarray}}where $C_{1}=3\Vert q\Vert _{\mathbf{V}}+\underset{m\in {\mathbb{N}_{0}}}\sup\Vert
p_{m}\Vert _{\mathbf{V}}+\left\vert \lambda \right\vert <\infty $, and $\hat{y}_{m}(x):=\underset{t\in \lbrack 0,x]}\max \left\vert y_{m}(t)\right\vert \in {\mathcal{C}}(I,{\mathbb{R}}%
)$. Obviously, $\hat{y}_{m}(x)$ is non-decreasing in $x\in I$. By substituting (%
\ref{zmwm}) and (\ref{wmym}) into (\ref{ymzm}), we have
\begin{equation*}
|y_{m}(t)|\leqslant C_{2}+C_{1}\int_{[0,x]}\hat{y}_{m}(s) \mathrm{d}%
s\qquad \forall  t\in
\lbrack 0,x]\subset I,
\end{equation*}%
where $C_{2}:=|y_{0}|+|z_{0}|+\frac{1}{2}|w_{0}|$. Thus{,}
\begin{equation*}
\hat{y}_{m}(x)=\max_{t\in \lbrack 0,x]}|y_{m}(t)|\leqslant
C_{2}+C_{1}\int_{[0,x]}\hat{y}_{m}(s) \mathrm{d}s,\qquad x\in I.
\end{equation*}%
{Then the Gronwall inequality together with the fact} $\hat{y}_{m}(0)=|y_{m}(0)|=|y_{0}|$ shows that $\underset{m\in {\mathbb{N}_{0}}}\sup\Vert \hat{y}_{m}\Vert _{\infty
}\leqslant C_{3}$, where {$C_{3}=C_{2}e^{C_{1}}$}%
. Hence, $\underset{m\in {\mathbb{N}_{0}}}\sup\Vert y_{m}\Vert _{\infty }\leqslant
\underset{m\in {\mathbb{N}_{0}}}\sup\Vert \hat{y}_{m}\Vert _{\infty }$
$\leqslant C_{3}$%
.

\textit{Step 2.} {Our task now is to prove} the sequence $\{y_{m}\}_{m\in {\mathbb{N}_{0}}}$ is relatively compact in $({\mathcal{C}}(I,{\mathbb{C}}),\Vert \cdot \Vert _{{\mathcal{\infty }}})$.

{The equation} (\ref{fix90}) {leads to}
\begin{eqnarray}
&&y_{m}(x)=\tilde{y}_{0}(x)-{\mathcal{Z}}\left( p_{m},y_{m}\right) (x),  \quad x \in I{,} \label{fix901} \\
&&{{\mathcal{Z}}\left( p_{m},y_{m}\right)}:=\int_{I}F(x,t)q\left(
t\right) y_{m}(t) \mathrm{d}t-\int_{I}G(x,t)y_{m}(t) \mathrm{d}\tilde{\mu}%
_{m}(t),x \in I, \qquad \label{zpym}
\end{eqnarray}
where $\tilde{\mu}_{m}(x)=q(x)+ip_{m}(x)-\lambda ix$. For any $0\leqslant x_{1}\leqslant x_{2}\leqslant 1$, {the following identity is obtained} from (\ref{fix901}),
\begin{eqnarray*}
|y_{m}(x_{2})-y_{m}(x_{1})|\!\!\!\! &\leqslant &\!\!\!\!\left( \left\vert z_{0}\right\vert +\left\vert
w_{0}\right\vert \right) \left\vert x_{2}-x_{1}\right\vert +2\left\vert
x_{2}-x_{1}\right\vert \Vert q\Vert _{\infty }\sup_{m\in {\mathbb{N}_{0}}%
}\Vert y_{m}\Vert _{\infty } \\
&&\!\!\!\!+\left\vert x_{2}-x_{1}\right\vert \underset{m\in {\mathbb{N}_{0}}}\sup\Vert
y_{m}\Vert _{\infty }\left( \Vert q\Vert _{\mathbf{V}}+\sup_{m\in {\mathbb{%
N}}}\Vert p_{m}\Vert _{\mathbf{V}}+\left\vert \lambda \right\vert \right)
\\
\!\!\!\! &\leqslant &\!\!\!\!C_{4}\left\vert x_{2}-x_{1}\right\vert ,
\end{eqnarray*}%
where $C_{4}:=\left\vert z_{0}\right\vert +\left\vert w_{0}\right\vert
+C_{1}C_{3}$. Hence{,} $\{y_{m}\}_{m\in {\mathbb{N}_{0}}}$ is {equicontinuous}. From {Arzel\`{a}-Ascoli} theorem, there exists a subsequence $\{y_{{m_{k}}}\}_{{k}\in {\mathbb{N}_{0}}}$ of $\{y_{m}\}_{m\in {\mathbb{N}_{0}}%
} $ such that $\{y_{{m_{k}}}\}_{{k}\in {\mathbb{N}_{0}}}$ converges uniformly to a continuous function $y_{\ast }$.

\textit{Step 3.} {We have to show that} the sequence $\{y_{m}\}_{m\in {\mathbb{N}_{0}}}$ is relatively compact in\/ $({\mathcal{C}}^{1}(I,{\mathbb{C}}),\Vert \cdot
\Vert _{{\mathcal{C}}^{1}})$.

For each ${k}\in {\mathbb{N}_{0}}$, $y_{{m_{k}}}$ is
continuously differentiable, and $y_{{m_{k}}}^{\prime }=z_{{m_{k}}}$%
. {With} $\vert z_{{m_{k}}}(0)\vert=\vert z_{0}\vert$, {substitution of} (\ref{wmym}) into (\ref{zmwm}) {yields}
\begin{equation*}
|z_{{m_{k}}}(x)|\leqslant C_{5},\qquad x\in I ,
\end{equation*}%
where $C_{5}:=|z_{0}|+|w_{0}|+C_{1}C_{3}$. Therefore, $\underset{{k}\in {\mathbb{N}_{0}}}\sup\Vert z_{{m_{k}}}\Vert _{\infty }\leqslant C_{5}$, {{i.e.}},
the sequence $\{y_{{m_{k}}}^{\prime }\}_{{k}\in {\mathbb{N}_{0}}}$
is uniformly bounded. {The following identity is obtained from (\ref{fix901}),}
\begin{eqnarray*}
y_{{m_{k}}}^{\prime }(x)\!\!\!\! &=&\!\!\!\!\tilde{y}_{0}^{\prime
}(x)-\left( {\mathcal{Z}}\left( p_{{m_{k}}},y_{{m_{k}}}\right)
\right) ^{\prime }(x) \\
&=&\!\!\!\!z_{0}+w_{0}x-\int_{I}F^{\prime }(x,t)q\left( t\right)
y_{{m_{k}}}(t) \mathrm{d}t+\int_{I}G^{\prime }(x,t)y_{{m_{k}}}(t) %
\mathrm{d}\tilde{\mu}_{{m_{k}}}(t)%
.
\end{eqnarray*}%
For any $0\leqslant x_{1}\leqslant x_{2}\leqslant 1$, one has
\begin{eqnarray*}
|y_{{m_{k}}}^{\prime }(x_{2})-y_{{m_{k}}}^{\prime }(x_{1})|\!\!\!\!\!
&\leqslant &\!\!\!\!\!{|w_{0}|\left\vert x_{2}-x_{1}\right\vert +\left\vert x_{2}-x_{1}\right\vert \sup_{{k}\in {\mathbb{N}_{0}}}\Vert y_{{m_{k}}}\Vert _{\infty }\left(\sup_{{k}\in {\mathbb{N}_{0}}}\Vert p_{{m_{k}}}\Vert _{\mathbf{%
V}}+ \Vert q\Vert _{\mathbf{V}%
}+\left\vert \lambda \right\vert \right) } \\
&\leqslant &\!\!\!\!\!C_{6}\left\vert x_{2}-x_{1}\right\vert ,
\end{eqnarray*}%
where $C_{6}:=|w_{0}|+C_{1}C_{3}$. Hence, $\{y_{{m_{k}}}^{\prime
}\}_{{k}\in {\mathbb{N}_{0}}}$ is {equicontinuous}. According to {Arzel\`{a}-Ascoli} theorem, there exists a subsequence $\{y_{{m_{k_{h}}}}^{\prime }\}_{{h}\in {\mathbb{N}_{0}}}$ of $\{y_{{m_{k}}}^{\prime }\}_{{k}\in {\mathbb{N}_{0}}}$ such that $\{y_{{m_{k_{h}}}}^{\prime }\}_{{h}\in {\mathbb{N}_{0}}}$ is uniformly
convergent to a continuous function $z_{\ast }$. Therefore, $y_{\ast }$ is
continuously differentiable, and
\begin{equation*}
y_{\ast }^{\prime }(x)=z_{\ast }(x),\qquad x\in I.
\end{equation*}%
{This implies that} the sequence $\{y_{m}\}_{m\in {\mathbb{N}_{0}}}$ is
relatively compact in $({\mathcal{C}}^{1}(I,{\mathbb{C}}),\Vert
\cdot \Vert _{{\mathcal{C}}^{1}})$. \qed

Now we turn to {prove} \PP \ref{wc}.

\smallskip\smallskip\noindent\emph{Proof of \PP \ref{wc}.} {$(i)$} {For any subsequence $\{y_{{m_{k}}}\}_{{k}\in {\mathbb{N}_{0}}}$ of $\{y_{m}\}_{m\in {\mathbb{N}_{0}}}$, it follows from Lemma \ref{rc} that} there is a sub-subsequence $\{y_{{m_{k_{h}}}}\}_{{h}\in {\mathbb{N}_{0}}}$ such that
\begin{equation}
y_{{m_{k_{h}}}}\rightarrow y_{\ast }\mbox{ in }({\mathcal{C}}^{1}(I,{\mathbb{C}}%
),\Vert \cdot \Vert _{{\mathcal{C}}^{1}})  \label{Y0}
\end{equation}%
for some $y_{\ast }\in {\mathcal{C}}^{1}(I,{\mathbb{C}})$. Let
\begin{equation*}
{\mathcal{Z}}\left( p_{0},y\right) \left( x\right) :=\int_{I}F(x,t)q\left( t\right) y(t) %
\mathrm{d}t-\int_{I}G(x,t)y(t) \mathrm{d}\tilde{\mu}_{0}(t),\quad  x\in I,
\end{equation*}
where $\tilde{\mu}_{0}(x)=q(x)+ip_{0}(x)-\lambda ix$. From (\ref{zpym}), one has
\begin{eqnarray*}
{\mathcal{Z}}\left( p_{{m_{k_{h}}}},y_{{m_{k_{h}}}}\right) (x)-{\mathcal{Z}}\left(
p_{0},y_{\ast }\right) (x)\hh\!\! &=&\hh\!\!\left[ \int_{I}F(x,t)q\left( t\right) \left(
y_{{m_{k_{h}}}}(t)-y_{\ast }(t) \right) \mathrm{d}t \right.\\
&&-\int_{I}G(x,t)\left(y_{{m_{k_{h}}}}(t) -y_{{\ast }}(t) \right) \mathrm{d}\tilde{\mu}_{{m_{k_{h}}}}(t)%
\biggr]\\
&&-i\left[ \int_{I}G(x,t) y_{\ast }(t) \mathrm{d}p_{{{m_{k_{h}}}}}(t)\right.
\\
&&\left. {-}\int_{I}G(x,t) y_{\ast }(t) \mathrm{d}p_{0}\left( t\right) %
\right] \\
&=:&\!\!\!\!J_{{m_{k_{h}}}}(x)+iK_{{m_{k_{h}}}}(x).
\end{eqnarray*}%
From (\ref{Y0}), it yields
\begin{eqnarray}
|J_{{m_{k_{h}}}}(x)|\!\!\!\! &\leqslant &\!\!\!\!2\Vert q\Vert _{\mathbf{V}}\Vert
y_{{m_{k_{h}}}}-y_{\ast }\Vert _{\infty }+\frac{1}{2}\Vert y_{{m_{k_{h}}}}-y_{\ast
}\Vert _{\infty }\left( \Vert q\Vert _{\mathbf{V}}+\Vert p_{{m_{k_{h}}}}\Vert _{%
\mathbf{V}}+\left\vert \lambda \right\vert \right) \notag \\
&\leqslant &\!\!\!\!\left( 2\Vert q\Vert _{\mathbf{V}}+\frac{1}{2}\left(
\Vert q\Vert _{\mathbf{V}}+\sup_{{l}\in {\mathbb{N}_{0}}}\Vert
p_{{m_{k_{h}}}}\Vert _{\mathbf{V}}+\left\vert \lambda \right\vert \right) \right)
\Vert y_{{m_{k_{h}}}}-y_{\ast }\Vert _{{\mathcal{C}}^{1}}\notag\\
&\rightarrow & 0 \qquad \text{ as }{h}\rightarrow \infty . \label{jm3}
\end{eqnarray}%
For any fixed $x \in I$, {$G(x,\cdot )y_{\ast }(\cdot )\in \mathcal{C}(I,\mathbb{C})$}. Since $p_{m}\rightarrow
p_{0}$ in $({\mathcal{M}}_{0}(I,{\mathbb{K}}),w^{\ast })$, for $x\in I$, one has $K_{{m_{k_{h}}}}(x)\rightarrow 0$ as $h\rightarrow 0$. Therefore,
\begin{equation}
\lim_{{h}\rightarrow \infty }{\mathcal{Z}}\left(
p_{{m_{k_{h}}}},y_{{m_{k_{h}}}}\right) (x)={\mathcal{Z}}\left( p_{0},y_{\ast }\right)
(x)\!\!\!\!\qquad \mbox{for each }x\in I.  \label{z01}
\end{equation}
From the equality (\ref{fix901}), the uniform convergence {in} (\ref{Y0}) and the {pointwise} convergence {in} (\ref{z01}), we have
\begin{equation*}
y_{\ast }(x)=\tilde{y}_{0}(x)-{\mathcal{Z}}\left( p_{0},y_{\ast }\right)
(x),\qquad x\in I.
\end{equation*}%
{Then it follows from (\ref{fix90}) that $y_{\ast }(x)=y_{0}(x)=y(x,\lambda ,p_{0},q)$. Since} the limit $y_{\ast }=y_{0}$ is independent of the choice of ${m_{k_{h}}}$, it yields that $%
y_{m}\rightarrow y_{0}$ in $({\mathcal{C}}^{1}(I,{\mathbb{C}}),\Vert \cdot
\Vert _{{\mathcal{C}}^{1}})$, {{i.e.}},
\begin{equation*}
\underset{m\rightarrow \infty }{\lim }\Vert y(x,\lambda
,p_{m},q)-y(x,\lambda ,p_{0},q)\Vert _{{\mathcal{C}}^{1}}=0;
\end{equation*}%
this proves the continuity {in} (\ref{py}), and
\begin{equation*}
\underset{m\rightarrow \infty }{\lim }\Vert y^{\prime }(x,\lambda
,p_{m},q)-y^{\prime }(x,\lambda ,p_{0},q)\Vert _{\infty }=0.
\end{equation*}

Next, for $f$, $g\in {\mathcal{C}}(I,{\mathbb{K}})$ and $F(x):=\int_{[0,x]}f(t) %
\mathrm{d}\tilde{\mu}(t)$, $x \in I$, from \cite[p. 260, Theorem G]{Mi83}, { we get} the equality
\begin{equation*}
\int_{(0,1]}g(t)\mathrm{d}F(t)=\int_{(0,1]}g(t)f(t) \mathrm{d}\tilde{\mu}%
(t).
\end{equation*}%
Then for $m\in {\mathbb{N}_{0}}$ and $f\in {%
\mathcal{C}}(I,{\mathbb{K}}),$ we obtain
\begin{eqnarray*}
\int_{I}f(x)\mathrm{d}\left( y_{m}^{\prime }\right) ^{\bullet
}(x)\!\!\!\! &=&\!\!\!\!f(0)(\left( y_{m}^{\prime }\right) ^{\bullet
}(0+)-\left( y_{m}^{\prime }\right) ^{\bullet }(0))+\int_{(0,1]}f(x)%
\mathrm{d}\left( y_{m}^{\prime }\right) ^{\bullet }(x) \\
&=&\!\!\!\!-f(0)y_{m}\left( 0\right) \mu
_{m}(0+)-\int_{(0,1]}f(x)2q(x)z_{m}(x) \mathrm{d}x \\
&&\hh -\int_{(0,1]}f(x) y_{m}(x) \mathrm{d}\mu _{m}(x) \\
&=&\!\!\!\!-\int_{I}f(x)2q(x)z_{m}(x) \mathrm{d}x-\int_{I}f(x)%
 y_{m}(x) \mathrm{d}\mu _{m}(x),
\end{eqnarray*}%
where $\mu _{m}(x)=q(x)-ip_{m}(x)+\lambda ix$, $x\in I$. When $m\rightarrow
\infty ,$ we obtain
\begin{eqnarray}
\int_{I}f(x)\mathrm{d}\left( y_{m}^{\prime }\right) ^{\bullet
}(x)\!\!\!\! &=&\!\!\!\!-\int_{I}f(x)2q(x)\left( z_{m}(x)-z_{0}\left( x\right)
\right)  \mathrm{d}x-\int_{I}f(x)2q(x)z_{0}\left( x\right) \mathrm{d}x
\notag \\
&&\!\!\!\!-\int_{I}f(x) \left( y_{m}(x)-y_{0}(x)\right)  \mathrm{d}\mu
_{m}(x)-\int_{I}f(x)y_{0}(x)  \mathrm{d}\mu _{m}(x)  \notag \\
&\rightarrow &\!\!\!\!-\int_{I}f(x)2q(x)z_{0}\left( x\right) \mathrm{d}%
x-\int_{I}f(x)y_{0}(x)  \mathrm{d}\mu _{0}(x)  \label{y21} \\
&=&\!\!\!\!\int_{I}f(x)\mathrm{d}\left( y_{0}^{\prime }\right)^
{\bullet }(x),  \notag
\end{eqnarray}%
{{i.e.}}, $\left( y_{m}^{\prime }\right) ^{\bullet }\rightarrow \left(
y_{0}^{\prime }\right) ^{\bullet }$ in $({\mathcal{M}}(I,{\mathbb{C}}%
),w^{\ast })$. This proves the continuity {in} (\ref{coy2}).

Let $f(x)\equiv 1$, then from (\ref{y21}) and
\begin{equation*}
\int_{I} \mathrm{d}\mu (x)=\mu (1)-\mu (0),
\end{equation*}
it yields $%
\underset{m\rightarrow \infty }\lim({\left( y_{m}^{\prime }\right)^{\bullet } (1)-\left(
y_{m}^{\prime }\right) ^{\bullet }(0)})=\left( y_{0}^{\prime }\right) ^{\bullet
}(1)-\left( y_{0}^{\prime }\right) ^{\bullet }(0)$. Since $\left(
y_{m}^{\prime }\right) ^{\bullet }(0)=\left( y_{0}^{\prime }\right)^
{\bullet }(0)=w_{0}$ {holds} for all $m\in\mathbb{N}$, we obtain that $\underset{m\rightarrow \infty }\lim\left( y_{m}^{\prime }\right) ^{\bullet }(1)=\left( y_{0}^{\prime
}\right) ^{\bullet }(1)$. This proves the continuity result {in}
(\ref{coy21}).

{$(ii)$ Suppose the sequence $\{q_{m}\}_{m \in {\mathbb{N}}}$ converges to $q_{0}$ in $({\mathcal{M}}_{0}(I,{\mathbb{K}}),w^{\ast})$. For $m \in {\mathbb{N}_{0}}$, let $y_{q_{m}}:=y(x,\lambda, p,q_{m})$, then following the same procedure as in the proof of Lemma \ref{rc}, we can prove that $\{y_{q_{m}}\}_{m \in {\mathbb{N}_{0}}}$ is relatively compact in the space $({\mathcal{C}}^{1}(I,{\mathbb{C}}),\Vert
\cdot \Vert _{{\mathcal{C}}^{1}})$. For any subsequence $\{y_{q_{{{m_{k}}}}}\}_{{k}\in {\mathbb{N}_{0}}}$ of $\{y_{q_{m}}\}_{m\in {\mathbb{N}_{0}}}$, we select a sub-subsequence $\{y_{q_{{m_{k_{h}}}}}\}_{{h}\in {\mathbb{N}_{0}}}$ such that
\begin{equation}
y_{q_{{m_{k_{h}}}}}\rightarrow y_{q\ast }\mbox{ in }({\mathcal{C}}^{1}(I,{\mathbb{C}}%
),\Vert \cdot \Vert _{{\mathcal{C}}^{1}}) \label{yqmqs}
\end{equation}%
for some $y_{q\ast }\in {\mathcal{C}}^{1}(I,{\mathbb{C}})$. Denote
\begin{equation*}
{\mathcal{Z}}\left( q_{m},y\right) \left( x\right) :=\int_{I}F(x,t)q_{m}\left( t\right) y(t)\mathrm{d}t-\int_{I}G(x,t)y(t) \mathrm{d}\tilde{\mu}_{q_{m}}(t),\quad m\in \mathbb{N}_{0},
\end{equation*}
where $\tilde{\mu}%
_{q_{m}}(x)=q_{m}(x)+ip(x)-\lambda ix$, $x\in I$. Then%
\begin{eqnarray*}
&&{\mathcal{Z}}\left( q_{{{m_{k_{h}}}}},y_{q_{{m_{k_{h}}}}}\right) (x)-{\mathcal{Z}}\left(
q_{0},y_{q\ast }\right) (x)\\
\!\!\!\! &=&\!\!\!\! \left[\int_{I}F(x,t)q_{{m_{k_{h}}}}( t) (y_{q_{{m_{k_{h}}}}}(t)-y_{q\ast }(t)) \mathrm{d}t +\int_{I}F(x,t)\left(q_{{m_{k_{h}}}}( t)-q_{0}( t) \right)y_{q\ast }(t) \mathrm{d}t\right.\\
&&\!\!\!\!-\int_{I}G(x,t)(y_{q_{{m_{k_{h}}}}}(t)-y_{q\ast }(t)) \mathrm{d}\tilde{\mu}_{q_{{m_{k_{h}}}}}(t)\biggr]\\
&&\!\!\!\!-\left[ \int_{I}G(x,t) y_{q\ast }(t) \mathrm{d}q_{{m_{k_{h}}}}(t)-\int_{I}G(x,t) y_{q\ast }(t) \mathrm{d}q_{0}\left( t\right) \right].
\end{eqnarray*}%
Here, {using the integration by parts formula} and the fact $q_{{m_{k_{h}}}}( 0)=q_{0}( 0)=0$, we have
\begin{eqnarray*}
&&\int_{I}F(x,t)\left(q_{{m_{k_{h}}}}( t)-q_{0}( t) \right)y_{q\ast }(t) \mathrm{d}t\\
\hh&=&\hh\int_{I}\frac{\mathrm{d}{\int_{[0,t]}F(x,s)y_{q\ast }(s)\mathrm{d}s}}{\mathrm{d}t}\left(q_{{m_{k_{h}}}}( t)-q_{0}( t) \right) \mathrm{d}t\\
\hh&=&\hh\left[\left(q_{{m_{k_{h}}}}( t)-q_{0}( t) \right)\int_{[0,t]}F(x,s)y_{q\ast }(s)\mathrm{d}s\right]\Bigg\vert_{t=0}^{t=1}-\int_{I}{\int_{[0,t]}F(x,s)y_{q\ast }(s)\mathrm{d}s}\mathrm{d}\left(q_{{m_{k_{h}}}}( t)-q_{0}( t) \right)\\
\hh&=&\hh\left(\int_{I}\mathrm{d}q_{{m_{k_{h}}}}( t)-\int_{I}\mathrm{d}q_{0}( t) \right)\int_{[0,1]}F(x,s)y_{q\ast }(s)\mathrm{d}s-\left[\int_{I}{\int_{[0,t]}F(x,s)y_{q\ast }(s)\mathrm{d}s}\mathrm{d}q_{{m_{k_{h}}}}( t)\right.\\
&&\hh-\left.\int_{I}{\int_{[0,t]}F(x,s)y_{q\ast }(s)\mathrm{d}s}\mathrm{d}q_{0}( t) \right].
\end{eqnarray*}
Note that for any fixed $x\in I$, $\int_{[0,t]}F(x,s)y_{q\ast }(s)\mathrm{d}s$ and $G(x,t) y_{q\ast }(t)$ are continuous functions {of} $t \in I$. Thus, {from} (\ref{yqmqs}) and {the fact} $q_{m}\rightarrow q_{0}$ in $({\mathcal{M}}_{0}(I,{\mathbb{K}}),w^{\ast })$, it yields
\begin{eqnarray}
\hh&&\hh\Big\vert {\mathcal{Z}}\left( q_{{{m_{k_{h}}}}},y_{q_{{m_{k_{h}}}}}\right) (x)-{\mathcal{Z}}\left(q_{0},y_{q\ast }\right) (x)\Big\vert  \notag\\
\hh&\leqslant&\hh2\sup_{{l}\in {\mathbb{N}_{0}}}\Vert q_{{m_{k_{h}}}}\Vert\Vert y_{q_{{m_{k_{h}}}}}-y_{q{\ast }}\Vert _{{\mathcal{C}}^{1}}+2\Vert y_{q^{\ast }}\Vert _{\infty}\Big\vert\int_{I}\mathrm{d}q_{{m_{k_{h}}}}( t)-\int_{I}\mathrm{d}q_{0}( t)\Big\vert\notag\\
\hh&&\hh+\Big\vert\int_{I}{\int_{[0,t]}F(x,s)y_{q\ast }(s)\mathrm{d}s}\mathrm{d}q_{{m_{k_{h}}}}( t)-\int_{I}{\int_{[0,t]}F(x,s)y_{q\ast }(s)\mathrm{d}s}\mathrm{d}q_{0}( t) \Big\vert \notag\\
\hh&&\hh+\frac{1}{2}\left(\sup_{{l}\in {\mathbb{N}_{0}}}\Vert q_{{m_{k_{h}}}}\Vert_{\mathbf{V}}+\Vert p\Vert_{\mathbf{V}}+\vert \lambda \vert\right)\Vert y_{q_{{m_{k_{h}}}}}-y_{q{\ast }}\Vert _{{\mathcal{C}}^{1}}\notag\\
\hh&&+\Big\vert \int_{I}G(x,t) y_{q\ast }(t) \mathrm{d}q_{{m_{k_{h}}}}(t) -\int_{I}G(x,t) y_{q\ast }(t) \mathrm{d}q_{0}\left( t\right)\Big\vert\notag\\
\hh&\rightarrow&\hh0 \qquad \text{as } h\rightarrow \infty \label{zqmzq0},
\end{eqnarray}
{{i.e.}},
\begin{equation*}
\lim_{{h}\rightarrow \infty }{\mathcal{Z}}\left( q_{{{m_{k_{h}}}}},y_{q_{{m_{k_{h}}}}}\right) (x)={\mathcal{Z}}\left(q_{0},y_{q\ast }\right) (x) \qquad \text{for each } x\in I.
\end{equation*}
Then \PP \ref{wc} $(ii)$ can be proved by an argument similar to the one used in \PP \ref{wc} $(i)$.
}\qed

\begin{rem}
\label{uc1}
{It should be mentioned that the continuity in $(\ref{py})$ and $(\ref{yq})$ hold uniformly for $\lambda \in U$, where $U$ is any bounded subset of $\mathbb{C}$. Let $C_{U}:=\underset{\lambda \in U}{max}\vert \lambda \vert$. Note that the proofs of Lemma $\ref{rc}$ and \PP \ref{wc} go through if we replace $\vert \lambda \vert$ in {the definition of }$C_{1}$, $(\ref{jm3})$ and $(\ref{zqmzq0})$ by $C_{U}$. This implies that the relatively compactness of the sequences $\{y_{m}\}_{m \in \mathbb{N}_{0}}$, $\{y_{q_{m}}\}_{m \in \mathbb{N}_{0}}$ hold uniformly on $U$, and then we acquire the uniform continuity in $(\ref{py})$ and $(\ref{yq})$ for $\lambda \in U$. }
\end{rem}
We now construct an example to illustrate the continuity result {in} (\ref{coy21}) cannot be generalized to other $x \in (0,1)$.
\begin{ex}
{Suppose $\lambda =0$, $q=0$ {and} $(y_{0},z_{0},w_{0}) =(1,0,0)$. For $m \in \mathbb{N}$, let
\begin{equation*}
p_{m}(x):=\left\{
\begin{array}{ll}
0 & \mbox{ for }x\in \lbrack 0,\frac{1}{2}), \\
m(x-\frac{1}{2}) & \mbox{ for }x\in \lbrack \frac{1}{2},\frac{1}{2}+\frac{1}{%
m}), \\
1 & \mbox{ for }x\in \lbrack \frac{1}{2}+\frac{1}{m},1],%
\end{array}%
\right.
\end{equation*}%
then $p_{m}\rightarrow p_{0}$ in $({\mathcal{M}}_{0}(I,{\mathbb{R}}),w^{\ast
})$, where
\begin{equation*}
p_{0}(x)=\delta _{\frac{1}{2}}(x)=\left\{
\begin{array}{l}
0\mbox{ for }x\in \lbrack 0,\frac{1}{2}), \\
1\mbox{ for }x\in \lbrack \frac{1}{2},1].%
\end{array}%
\right.
\end{equation*}%
A simple calculation gives
\begin{equation*}
\lim_{m\rightarrow \infty }\left( y^{\prime }\right) ^{\bullet
}(1,0,p_{m},0)=i=\left( y^{\prime }\right) ^{\bullet }(1,0,\delta _{1/2},0),
\end{equation*}%
but
\begin{equation*}
\lim_{m\rightarrow \infty }\left( y^{\prime }\right) ^{\bullet }(\frac{1}{2}%
,0,p_{m},0)=0\neq i=\left( y^{\prime }\right) ^{\bullet }(\frac{1}{2}%
,0,\delta _{\frac{1}{2}},0).
\end{equation*}
}
\end{ex}

{
\begin{prop}
\label{uc2}
{$(i)$ Let $U$ be any bounded subset of $\mathbb{C}$, then the following mappings are uniformly continuous for $\lambda \in U$,}
\begin{eqnarray*}
&&({\mathcal{M}}_{0}(I,{\mathbb{K}}),\Vert \cdot \Vert _{\infty})\times ({\mathcal{M}}_{0}(I,{\mathbb{K}}),\Vert \cdot \Vert _{\infty})\rightarrow \left( C^{1}(I,{\mathbb{C}}),\Vert \cdot \Vert _{{\mathcal{C}}^{1}}\right) ,\quad
(p,q)\rightarrow y(\cdot ,\lambda ,p,q),\\
&&({\mathcal{M}}_{0}(I,{\mathbb{K}}),\Vert \cdot \Vert _{\infty})\times ({\mathcal{M}}_{0}(I,{\mathbb{K}}),\Vert \cdot \Vert _{\infty})\rightarrow \left( {\mathcal{M}}(I,{\mathbb{C}}),w^{\ast }\right) ,\quad\quad\ \
(p,q)\rightarrow (y^{\prime})^{\bullet}(\cdot ,\lambda ,p,q).
\end{eqnarray*}
{More precisely, for any $p_{0}$, $q_{0} \in ({\mathcal{M}}_{0}(I,{\mathbb{K}}),\Vert \cdot \Vert _{\infty})$ and $\epsilon >0$, there is a $\delta >0$ such that if $\Vert p-p_{0} \Vert_{\infty}+\Vert q-q_{0} \Vert_{\infty}\leqslant \delta$, one has
\begin{eqnarray}
&&\vert y(x,\lambda ,p,q)-y(x,\lambda ,p_{0},q_{0})\vert <\epsilon,\label{y0uc}\\
 &&\vert y^{\prime}(x,\lambda ,p,q)-y^{\prime}(x,\lambda ,p_{0},q_{0})\vert <\epsilon,\label{y1uc}\\
 &&\vert (y^{\prime})^{\bullet}(x,\lambda ,p,q)-(y^{\prime})^{\bullet}(x,\lambda ,p_{0},q_{0})\vert <\epsilon \label{y2uc}
\end{eqnarray}
hold uniformly for $x \in I$ and $\lambda \in U$.}

{$(ii)$ The following mappings are uniformly continuous for $\lambda \in U$,}
\begin{eqnarray*}
&&({\mathcal{M}}_{0}(I,{\mathbb{K}}),\Vert \cdot \Vert _{\mathbf{V}})\times ({\mathcal{M}}_{0}(I,{\mathbb{K}}),\Vert \cdot \Vert _{\mathbf{V}})\rightarrow \left( C^{1}(I,{\mathbb{C}}),\Vert \cdot \Vert _{{\mathcal{C}}^{1}}\right) ,\quad
(p,q)\rightarrow y(\cdot ,\lambda ,p,q),\\
&&({\mathcal{M}}_{0}(I,{\mathbb{K}}),\Vert \cdot \Vert _{\mathbf{V}})\times ({\mathcal{M}}_{0}(I,{\mathbb{K}}),\Vert \cdot \Vert _{\mathbf{V}})\rightarrow \left( {\mathcal{M}}(I,{\mathbb{C}}),w^{\ast }\right) ,\quad\quad\ \
(p,q)\rightarrow (y^{\prime})^{\bullet}(\cdot ,\lambda ,p,q).
\end{eqnarray*}
{That is to say, for any $p_{0}$, $q_{0} \in ({\mathcal{M}}_{0}(I,{\mathbb{K}}),\Vert \cdot \Vert _{\mathbf{V}})$ and $\epsilon >0$, there is a $\delta >0$ such that if $\Vert p-p_{0} \Vert_{\mathbf{V}}+\Vert q-q_{0} \Vert_{\mathbf{V}}\leqslant \delta$, the inequalities $($\ref{y0uc}$)$-$($\ref{y2uc}$)$ hold uniformly for $x \in I$ and $\lambda \in U$.}
\end{prop}
\Proof $(i)$ Suppose the sequence $\{p_{m}\}_{m\in \mathbb{N}}$ converges to $p_{0}$ in $({\mathcal{M}}_{0}(I,{\mathbb{K}}),\Vert \cdot \Vert _{\infty})$, and the sequence $\{q_{m}\}_{m\in \mathbb{N}}$ converges to $q_{0}$ in $({\mathcal{M}}_{0}(I,{\mathbb{K}}),\Vert \cdot \Vert _{\infty})$, then there are constants $M_{p}$, $M_{q}$ such that
\begin{equation}
\sup_{m\in {\mathbb{N}}}\Vert p_{m}\Vert_{\mathbf{V}}<M_{p}, \qquad \sup_{m\in {\mathbb{N}}}\Vert q_{m}\Vert_{\mathbf{V}}<M_{q}.\label{Mpq}
\end{equation}
Due to {the} equations (\ref{ryz})-(\ref{rwy}), $z=y^{\prime}$, and $w=(y^{\prime})^{\bullet}$, we have
\begin{eqnarray*}
\left(
\begin{array}{c}
y(x)\\
y^{\prime}(x)\\
(y^{\prime})^{\bullet}(x)
\end{array}
\right)\hh&=&\hh\left(
\begin{array}{c}
y_{0}\\
z_{0}\\
w_{0}
\end{array}
\right)+\left(
\begin{array}{c}
\int_{[0,x]}y^{\prime}(t)\mathrm{d}t\\
\int_{[0,x]}(y^{\prime})^{\bullet}(t)\mathrm{d}t\\
-\int_{[0,x]}2q(t)y^{\prime}(t)\mathrm{d}t-\int_{[0,x]}y(t)\mathrm{d}\mu(t)
\end{array}
\right)\\
\hh&:=&\hh\left(
\begin{array}{c}
y_{0}\\
z_{0}\\
w_{0}
\end{array}
\right)+\int_{[0,x]}\mathrm{d}\left(
\begin{array}{ccc}
0 & t & 0 \\
0 & 0 & t \\
-\mu(t) &-2\int_{[0,t]}q(s)\mathrm{d}s&0
\end{array}%
\right)\left(
\begin{array}{c}
y(t)\\
y^{\prime}(t)\\
(y^{\prime})^{\bullet}(t)
\end{array}
\right).
\end{eqnarray*}
Let
\begin{equation*}
A_{m}(t):=\left(
\begin{array}{ccc}
0 & t & 0 \\
0 & 0 & t \\
-\mu_{p_{m},q_{m}}(t) &-2\int_{[0,t]}q_{m}(s)\mathrm{d}s&0
\end{array}%
\right),\quad \mu_{p_{m},q_{m}}(t)=q_{m}(t)-ip_{m}(t)+\lambda it, \quad m\in \mathbb{N}_{0},
\end{equation*}
and
\begin{equation*}
\vert A_{m}(t)\vert := \max \left\{\vert \mu_{p_{m},q_{m}}(t)\vert,\vert t\vert+\left\vert 2\int_{[0,t]}q_{m}(s)\mathrm{d}s\right\vert \right\}.
\end{equation*}
According to (\ref{Mpq}), there exists a constant $M_{A}$ such that
\begin{equation*}
\sup_{m\in {\mathbb{N}}}\Vert A_{m}\Vert_{\mathbf{V}}<M_{A}
\end{equation*}
holds uniformly for $\lambda \in U$. Then from the proof of \cite[Theorem 4.1]{HT09}, we can prove $(i)$.

$(ii)$ Using the fact that $\Vert f\Vert _{\infty }\leqslant \Vert f\Vert _{\mathbf{V}}$ for all $f\in {\mathcal{M}}_{0}(I,{\mathbb{K}})$, we can prove the statement in $(ii)$. \qed}

Note that \PP \ref{wc} illustrates the dependence of $y(x,\lambda ,p,q)$, $y^{\prime }(x,\lambda ,p,q)$ and $\left( y^{\prime }\right) ^{\bullet}(x,\lambda ,p,q)$ on $p$, ${q} \in ({\mathcal{M}}_{0}(I,{\mathbb{K}}),w^{\ast })$. Next, we prove that $y(x,\lambda ,p,q)$, $y^{\prime }(x,\lambda ,p,q)$ {and} $\left( y^{\prime }\right) ^{\bullet}(x,\lambda ,p,q)$ are continuous Fr\'{e}chet differentiable in $p$, $q\in ({\mathcal{M}}_{0}(I,{\mathbb{K}}),\Vert \cdot \Vert _{\mathbf{V}})$, respectively. {And} then we deduce the Fr\'{e}chet derivatives correspondingly. {We first introduce the definition of Fr\'{e}chet derivative and some notations which will be used in \PP \ref{main-dc}.
\begin{defn}
A  map $T$ from a Banach space $(W,\Vert \cdot \Vert_{W})$ into a Banach space $(Z, \Vert \cdot \Vert_{Z})$, $T:W\rightarrow Z$, is differentiable at a point $w\in W$ if there exists a bounded linear map $d_{w}T:W\rightarrow Z$ such that
\begin{equation*}
\Vert T(w+h)-T(w)-d_{w}T\cdot h\Vert_{W}=o(\Vert h \Vert_{Z}) \quad \text{as }h\rightarrow 0\text{ in }W.
\end{equation*}
Here, the map $d_{w}T$ is called the Fr\'{e}chet derivative of $T$ at $w$.
\end{defn}
For $\nu_{p}$, $\nu_{q}\in {\mathcal{M}}_{0}(I,{\mathbb{K}})$, let
\begin{eqnarray*}
&&\partial _{p}N_{p,q}(x)\cdot \nu _{p}:=\left(
\begin{array}{ccc}
\partial _{p}y_{1}(x,\lambda ,p,q)\cdot \nu _{p} & \partial _{p}y_{2}(x,\lambda ,p,q)\cdot \nu _{p} & \partial _{p}y_{3}(x,\lambda ,p,q)\cdot \nu _{p} \\
\partial _{p}y_{1}^{\prime }(x,\lambda ,p,q)\cdot \nu _{p} & \partial _{p}y_{2}^{\prime }(x,\lambda ,p,q)\cdot \nu _{p} &
\partial _{p}y_{3}^{\prime }(x,\lambda ,p,q)\cdot \nu _{p} \\
\partial _{p}\left( y_{1}^{\prime }\right) ^{\bullet }(x,\lambda ,p,q)\cdot \nu _{p} & \partial _{p}\left(
y_{2}^{\prime }\right) ^{\bullet }(x,\lambda ,p,q)\cdot \nu _{p} & \partial _{p}\left( y_{3}^{\prime
}\right) ^{\bullet }(x,\lambda ,p,q)\cdot \nu _{p}%
\end{array}%
\right), \\
&&\partial _{q}N_{p,q}(x)\cdot \nu _{q}:=\left(
\begin{array}{ccc}
\partial _{q}y_{1}(x,\lambda ,p,q)\cdot \nu _{q} & \partial _{q}y_{2}(x,\lambda ,p,q)\cdot \nu _{q} & \partial _{q}y_{3}(x,\lambda ,p,q)\cdot \nu _{q} \\
\partial _{q}y_{1}^{\prime }(x,\lambda ,p,q)\cdot \nu _{q} & \partial _{q}y_{2}^{\prime }(x,\lambda ,p,q)\cdot \nu _{q} &
\partial _{q}y_{3}^{\prime }(x,\lambda ,p,q)\cdot \nu _{q} \\
\partial _{q}\left( y_{1}^{\prime }\right) ^{\bullet }(x,\lambda ,p,q)\cdot \nu _{q} & \partial _{q}\left(
y_{2}^{\prime }\right) ^{\bullet }(x,\lambda ,p,q)\cdot \nu _{q} & \partial _{q}\left( y_{3}^{\prime
}\right) ^{\bullet }(x,\lambda ,p,q)\cdot \nu _{q}%
\end{array}%
\right).
\end{eqnarray*}
Here, for $i=1,2$, $j=1,2,3$, $\partial _{p}y_{j}^{(i-1)}(x,\lambda ,p,q)$ and $\partial _{p}\left( y_{j}^{\prime}\right) ^{\bullet }(x,\lambda ,p,q)$ denote the Fr\'{e}chet derivatives of $y_{j}^{(i-1)}(x,\lambda ,p,q)$, $\left( y_{j}^{\prime}\right) ^{\bullet }(x,\lambda ,p,q)\in ({\mathcal{M}}(I,{\mathbb{K}}),\Vert \cdot \Vert _{\mathbf{V}})$ at $p\in ({\mathcal{M}}_{0}(I,{\mathbb{K}}),\Vert \cdot \Vert _{\mathbf{V}})$, respectively. Similarly, for $i=1,2$, $j=1,2,3$, $\partial _{q}y_{j}^{(i-1)}(x,\lambda ,p,q)$ and $\partial _{q}\left( y_{j}^{\prime}\right) ^{\bullet }(x,\lambda ,p,q)$ denote the Fr\'{e}chet derivatives of $y_{j}^{(i-1)}(x,\lambda ,p,q)$, $\left( y_{j}^{\prime}\right) ^{\bullet }(x,\lambda ,p,q)\in ({\mathcal{M}}(I,{\mathbb{K}}),\Vert \cdot \Vert _{\mathbf{V}})$ at $q\in ({\mathcal{M}}_{0}(I,{\mathbb{K}}),\Vert \cdot \Vert _{\mathbf{V}})$, respectively.
}

\begin{prop}
\label{main-dc} {$(i)$ Let $x\in I$, $\lambda
\in {\mathbb{C}}$, $q\in {\mathcal{M}}_{0}(I,{\mathbb{K}})$ and $(y_{0},z_{0},w_{0})\in \mathbb{K}^{3}$ be fixed. Then $y(x,\lambda ,p,q)$, $y^{\prime }(x,\lambda
,p,q)$ and $\left( y^{\prime }\right) ^{\bullet }(x,\lambda,p,q) $ are continuously Fr\'{e}chet differentiable in $p\in ({\mathcal{M}}_{0}(I,{\mathbb{K}}),\Vert \cdot \Vert _{\mathbf{V}})$. {Moreover, for $x\in (0,1]$, $\nu
_{p}\in {\mathcal{M}}_{0}(I,{\mathbb{K}})$,}
\begin{equation}
\partial _{p}N_{p,q}(x)\cdot \nu
_{p}=iN_{p,q}(x)\int_{[0,x]}N_{p,q}^{-1}(t)\left(
\begin{array}{ccc}
0 & 0 & 0 \\
0 & 0 & 0 \\
y_{1}(t,\lambda ,p,q) & y_{2}(t,\lambda ,p,q) & y_{3}(t,\lambda ,p,q)%
\end{array}%
\right)  \mathrm{d}\nu _{p}(t). \label{dp}
\end{equation}}

$(ii)$ Let $x\in I$, $\lambda \in {\mathbb{C}}$, $p\in {\mathcal{M}}_{0}(I,{\mathbb{K}})$ and $(y_{0},z_{0},w_{0})\in \mathbb{K}^{3}$ be fixed. Then $y(x,\lambda ,p,q)$, $y^{\prime }(x,\lambda ,p,q)$ and $\left( y^{\prime }\right)^{\bullet }(x,\lambda ,p,q)$ are
continuously Fr\'{e}chet differentiable in $q\in ({\mathcal{M}}_{0}(I,{\mathbb{K}}),\Vert \cdot \Vert _{\mathbf{V}})$. {Moreover, for $x\in (0,1]$, $\nu
_{q}\in {\mathcal{M}}_{0}(I,{\mathbb{K}})$, }
\begin{eqnarray}
\partial _{q}N_{p,q}(x)\cdot \nu _{q}=&&\hh\hh -N_{p,q}(x)\left[
\int_{[0,x]}N_{p,q}^{-1}(t)\left(
\begin{array}{ccc}
0 & 0 & 0 \\
0 & 0 & 0 \\
y_{1}(t,\lambda ,p,q) & y_{2}(t,\lambda ,p,q) & y_{3}(t,\lambda ,p,q)%
\end{array}%
\right)  \mathrm{d}\nu _{q}(t)\right.  \notag \\
&&\hh\hh \left. +\int_{[0,x]}N_{p,q}^{-1}(t)\left(
\begin{array}{ccc}
0 & 0 & 0 \\
0 & 0 & 0 \\
y_{1}^{\prime }(t,\lambda ,p,q) & y_{2}^{\prime }(t,\lambda ,p,q) &
y_{3}^{\prime }(t,\lambda ,p,q)%
\end{array}%
\right)  2\nu _{q}(t)\mathrm{d}t
\right] .  \notag\\
&& \label{dq}
\end{eqnarray}
\end{prop}

\Proof $(i)$ {Denote $y_{1}(x):=y_{1}(x,\lambda ,p,q)$ and $z_{1}(x):=y_{1}(x,\lambda ,p+\nu_{p},q)$, then from (\ref{equation}), we have
\begin{eqnarray*}
&&i\mathrm{d}\left( y_{1}^{\prime }\right) ^{\bullet }+2iq\left(x\right) y_{1}^{\prime }\mathrm{d}x+y_{1}\left( i\mathrm{d}q\left( x\right) +\mathrm{d}p\left(x\right) \right) = \lambda y_{1}\mathrm{d}x, y_{1}(0)=1,y_{1}^{\prime }(0)=\left( y_{1}^{\prime }\right) ^{\bullet }(0)=0,\\
&&i\mathrm{d}\left( z_{1}^{\prime }\right) ^{\bullet }+2iq\left(x\right) z_{1}^{\prime }\mathrm{d}x+z_{1}\left( i\mathrm{d}q\left( x\right) +\mathrm{d}(p+\nu_{p})\left(x\right) \right) = \lambda z_{1}\mathrm{d}x, z_{1}(0)=1,z_{1}^{\prime }(0)=\left( z_{1}^{\prime }\right) ^{\bullet }(0)=0.
\end{eqnarray*}
Let $w_{1}(x):=z_{1}(x)-y_{1}(x)$, then $w_{1}(x)$ satisfies
\begin{equation*}
i\mathrm{d}\left( w_{1}^{\prime }\right) ^{\bullet }+2iq\left(x\right) w_{1}^{\prime }\mathrm{d}x+w_{1}\left( i\mathrm{d}q\left( x\right) +\mathrm{d}p\left(x\right)-\lambda \mathrm{d}x \right) =-z_{1}\mathrm{d}\nu_{p}(x), w_{1}(0)=w_{1}^{\prime }(0)=\left( w_{1}^{\prime }\right) ^{\bullet }(0)=0.
\end{equation*}
{From} Lemma \ref{voc0}, it yields
\begin{equation*}
\left(
\begin{array}{c}
w_{1}(x) \\
w_{1}^{\prime }(x) \\
\left( w_{1}^{\prime }\right) ^{\bullet }(x)%
\end{array}%
\right) =N_{p,q}(x)\int_{[0,x]}N_{p,q}^{-1}(t)\left(
\begin{array}{c}
0 \\
0 \\
iz_{1}(t)
\end{array}
\right)  \mathrm{d}\nu_{p}(t), \quad x\in (0,1].
\end{equation*}
Then according to {\PP}\ref{uc2} $(ii)$ and the fact that $y_{j}(x,\lambda ,p,q)$, $y^{\prime}_{j}(x,\lambda ,p,q)$, $\left( y_{j}^{\prime}\right) ^{\bullet }(x,\lambda ,p,q)$, $j=1,2,3$ are bounded on $I$, we can obtain
\begin{eqnarray*}
&&\hh\left(
\begin{array}{c}
z_{1}(x)-y_{1}(x) \\
z_{1}^{\prime }(x)-y_{1}^{\prime }(x) \\
\left( z_{1}^{\prime }\right) ^{\bullet }(x)-\left( y_{1}^{\prime }\right) ^{\bullet }(x)
\end{array}%
\right) -N_{p,q}(x)\int_{[0,x]}N_{p,q}^{-1}(t)\left(
\begin{array}{c}
0 \\
0 \\
iy_{1}(t)
\end{array}
\right)  \mathrm{d}\nu_{p}(t)\notag\\
\hh&=&\hh N_{p,q}(x)\int_{[0,x]}N_{p,q}^{-1}(t)\left(
\begin{array}{c}
0 \\
0 \\
i(z_{1}(t)-y_{1}(t))
\end{array}
\right)  \mathrm{d}\nu_{p}(t)\notag\\
\hh&=&\hh\left(
\begin{array}{c}
o(\Vert \nu_{p}\Vert _{\mathbf{V}}) \\
o(\Vert \nu_{p}\Vert _{\mathbf{V}})  \\
o(\Vert \nu_{p}\Vert _{\mathbf{V}})
\end{array}
\right) \quad \text{as }\nu_{p}\rightarrow 0 \text{ in } ({\mathcal{M}}_{0}(I,{\mathbb{K}}),\Vert \cdot \Vert _{\mathbf{V}}). \label{yvq}
\end{eqnarray*}
Thus, the differentiability of $y_{1}(x,\lambda ,p,q)$, $y_{1}^{\prime}(x,\lambda ,p,q)$ and $\left( y_{1}^{\prime}\right) ^{\bullet }(x,\lambda ,p,q)$ in $p$ can be proved, and their derivatives are also obtained. Similarly, we can prove that $y_{j}^{(i-1)}(x,\lambda ,p,q)$, $\left( y_{j}^{\prime}\right) ^{\bullet }(x,\lambda ,p,q)$, $i=1,2$, $j=2,3$ are differentiable in $p$, and thus the equality $(\ref{dp})$ is obtained.
}

$(ii)$ Proceeding as in the proof of {\PP \ref{main-dc}} $(i)$, we can prove {\PP \ref{main-dc}} $(ii)$. \qed

We remark that for the derivatives of solutions of ordinary differential equation{s} with respect to
{the coefficients}, formulas like (\ref{dp})-(\ref{dq}) can be found in \cite{amour2001siam, ptrbookisp,
MZ96, Ze05}.

\subsection{The Asymptotic Formulae and Analyticity of Solutions}

Now we deduce the estimates of solutions and the analytic dependence of solutions on the spectral parameter $\lambda$ when $p$, $q\in {\mathcal{M}}_{0}(I,{\mathbb{R}})$, $\lambda \in \mathbb{C}$. Recall the definition in Lemma \ref{Vmu1} and denote $\mathbf{\check{p}}\left( x\right) :=\int_{[0,x]}\mathrm{d}\mathbf{V}%
_{p}\left( t\right) $, $\mathbf{\check{q}}\left( x\right) :=\int_{[0,x]}%
\mathrm{d}\mathbf{V}_{q}\left( t\right) $, then $\mathbf{\check{p}}\left( 1\right)=\Vert p \Vert_{\mathbf{V}} $, $\mathbf{\check{q}}\left( 1\right)=\Vert q \Vert_{\mathbf{V}} $.

\begin{thm}
\label{eopq} {Let
\begin{equation*}
\Xi \left( x,\lambda \right) :=e^{\left(\left\vert \mathrm{Im}\frac{k}{2}\right\vert +\left\vert \mathrm{Im}\frac{\omega k}{2}\right\vert +\left\vert \mathrm{Im}\frac{\omega ^{2}}{2}k\right\vert \right) x}, \quad k=\lambda ^{\frac{1}{3}},\quad {\omega =e^{\frac{2}{3}\pi i}}.
\end{equation*}
For $(x,\lambda ,p,q)\in I \times {%
\mathbb{C}\times \mathcal{M}}_{0}(I,{\mathbb{R}})\times {\mathcal{M}}_{0}(I,{%
\mathbb{R}})$, $j=1, 2, 3$, {we have}
\begin{eqnarray}
&&\left\vert y_{j}(x,\lambda ,p,q)\right\vert \leqslant \frac{3}{\left\vert
k\right\vert ^{j-1}}\Xi \left( x,\lambda \right) e^{3\left( 2\Vert q \Vert_{\mathbf{V}}  +\mathbf{\check{p}}\left( x\right) +\mathbf{\check{q}}\left(
x\right) \right) },  \label{ypq}   \\
&&\left\vert y_{j}(x,\lambda ,p,q)-y_{j}(x,\lambda ,0,0)\right\vert \leqslant
\frac{3}{\left\vert k\right\vert ^{j}}\Xi \left( x,\lambda \right)
e^{3\left( 2\Vert q \Vert_{\mathbf{V}}  +\mathbf{\check{p}}\left( x\right) +\mathbf{%
\check{q}}\left( x\right) \right) }. \label{pq-00}
\end{eqnarray}
}
\end{thm}

Note that, when $(p,q)=(0,0)$, the equation (\ref{equation}) reduces to
\begin{equation}
i\mathrm{d}\left( y^{\prime }\right) ^{\bullet }-\lambda y\mathrm{d}x=0.  \label{00}
\end{equation}
In order to prove Theorem \ref{eopq}, we need some properties of the solutions of {(\ref{00})}, which can be found in \cite[Lemma 2.1-2.3]{amour1999siam}.

The fundamental solutions of $(\ref{00})$ are
\begin{eqnarray*}
y_{1}\left( x,\lambda ,0,0\right) \!\!\!\! &=&\!\!\!\!\frac{1}{3}\left(
e^{ikx}+e^{i\omega kx}+e^{i\omega ^{2}kx}\right) \\
&=&\!\!\!\!\frac{1}{3}\left( 4\cos \left( \frac{kx}{2}\right) \cos \left(
\frac{\omega kx}{2}\right) \cos \left( \frac{\omega ^{2}kx}{2}\right)
-1\right. \\
&&\left. -4i\sin \left( \frac{kx}{2}\right) \sin \left( \frac{\omega kx}{2}%
\right) \sin \left( \frac{\omega ^{2}kx}{2}\right) \right) ,
\end{eqnarray*}
\begin{eqnarray*}
y_{2}\left( x,\lambda ,0,0\right) \!\!\!\! &=&\!\!\!\!\frac{1}{3ki}\left(
e^{ikx}+\omega ^{2}e^{i\omega kx}+\omega e^{i\omega ^{2}kx}\right) \\
&=&\!\!\!\!\frac{1}{3ki}\left( 4\cos \left( \frac{kx}{2}\right) \cos \left(
\frac{\omega kx}{2}-\frac{\pi }{3}\right) \cos \left( \frac{\omega ^{2}kx}{2}%
+\frac{\pi }{3}\right) -1\right. \\
&&\left. -4i\sin \left( \frac{kx}{2}\right) \sin \left( \frac{\omega kx}{2}-%
\frac{\pi }{3}\right) \sin \left( \frac{\omega ^{2}kx}{2}+\frac{\pi }{3}%
\right) \right) ,\\
y_{3}\left( x,\lambda ,0,0\right) \!\!\!\! &=&\!\!\!\!\frac{1}{3\left(
ki\right) ^{2}}\left( e^{ikx}+\omega e^{i\omega kx}+\omega ^{2}e^{i\omega
^{2}kx}\right) \\
&=&\!\!\!\!\frac{1}{3\left( ki\right) ^{2}}\left( 4\cos \left( \frac{kx}{2}%
\right) \cos \left( \frac{\omega kx}{2}+\frac{\pi }{3}\right) \cos \left(
\frac{\omega ^{2}kx}{2}-\frac{\pi }{3}\right) -1\right. \\
&&\left. -4i\sin \left( \frac{kx}{2}\right) \sin \left( \frac{\omega kx}{2}+%
\frac{\pi }{3}\right) \sin \left( \frac{\omega ^{2}kx}{2}-\frac{\pi }{3}%
\right) \right) .
\end{eqnarray*}%
For $j=1,2,3,$ we have
\begin{equation}
\left\vert y_{j}(x,\lambda ,0,0)\right\vert \leqslant \frac{3}{\left\vert
k\right\vert ^{j-1}}\Xi \left( x,\lambda \right) . \label{222}
\end{equation}
According to the identities
\begin{eqnarray*}
y_{j}(x,\lambda ,0,0)\!\!\!\! &=&\!\!\!\!\int_{[0,x]}y_{j-1}(t,\lambda ,0,0)%
\mathrm{d}t, \\
y_{j}^{\prime }(x,\lambda ,0,0)\!\!\!\! &=&\!\!\!\!y_{j-1}(x,\lambda ,0,0),\
j=2,3,
\end{eqnarray*}%
we acquire
\begin{equation}
\left\vert \frac{\partial ^{m}}{\partial x^{m}}y_{j}(x,\lambda
,0,0)\right\vert \leqslant 3\Xi \left( x,\lambda \right), j=1,2,3, m\leqslant j-1, m\in {\mathbb{N}}. \label{3}
\end{equation}

\smallskip\smallskip\noindent\emph{Proof of Theorem \ref{eopq}.} {Recall the definition of $N_{p,q}(x)$, then we have}
\begin{equation*}
N_{0,0}(x):=\left(
\begin{array}{ccc}
y_{1}(x,\lambda ,0,0) & y_{2}(x,\lambda ,0,0) & y_{3}(x,\lambda ,0,0) \\
y_{1}^{\prime }(x,\lambda ,0,0) & y_{2}^{\prime }(x,\lambda ,0,0) &
y_{3}^{\prime }(x,\lambda ,0,0) \\
\left( y_{1}^{\prime }\right) ^{\bullet }(x,\lambda ,0,0) & \left(
y_{2}^{\prime }\right) ^{\bullet }(x,\lambda ,0,0) & \left( y_{3}^{\prime
}\right) ^{\bullet }(x,\lambda ,0,0)%
\end{array}%
\right) , x\in I,
\end{equation*}%
and $N_{0,0}(x)N_{0,0}^{-1}(t)=N_{0,0}(x-t)$.

Let us {re}write the differential equation (\ref{equation}) as an inhomogeneous
differential equation%
\begin{equation*}
 i\mathrm{d}\left( y^{\prime }\right) ^{\bullet }-\lambda y\mathrm{d}x=-2iq\left(
x\right) y^{\prime }\mathrm{d}x-y\left( i\mathrm{d}q\left( x\right) +\mathrm{d}p\left(
x\right) \right) .
\end{equation*}%
For $x\in \left( 0,1\right] $, by Lemma \ref{voc0}, the fundamental
solutions $y_{j}\left( x,\lambda ,p,q\right) $, $j=1, 2, 3$ satisfy the
{following formula,}
\begin{eqnarray}
\left(
\begin{array}{c}
y_{j}(x,\lambda ,p,q) \\
y_{j}^{\prime }(x,\lambda ,p,q) \\
\left( y_{j}^{\prime }\right) ^{\bullet }(x,\lambda ,p,q)%
\end{array}%
\right) \!\!\!\! &=&\!\!\!\!\left(
\begin{array}{c}
y_{j}(x,\lambda ,0,0) \\
y_{j}^{\prime }(x,\lambda ,0,0) \\
\left( y_{j}^{\prime }\right) ^{\bullet }(x,\lambda ,0,0)%
\end{array}%
\right)-\int_{[0,x]}N_{0,0}(x-t)\left(
\begin{array}{c}
0 \\
0 \\
y_{j}\left( t,\lambda ,p,q\right)%
\end{array}%
\right) \mathrm{d}\tilde{\mu}_{p,q}\left( t\right)\notag \\
&&  -\int_{[0,x]}N_{0,0}(x-t)\left(
\begin{array}{c}
0 \\
0 \\
2q\left( t\right) y_{j}^{\prime }\left( t,\lambda ,p,q\right)%
\end{array}%
\right) \mathrm{d}t  ,  \label{voc1}
\end{eqnarray}%
where $\tilde{\mu}_{p,q}\left( t\right) =q\left( t\right) -ip\left( t\right)
$. Since $y_{3}\left( x-t,\lambda ,0,0\right) =0$ for $x=t$, the formula (\ref{voc1}) is also true for $x=0$. {From (\ref{voc1}), we see
\begin{eqnarray*}
y_{j}(x,\lambda,p,q)\!\!\!\!&=&\!\!\!\!y_{j}(x,\lambda, 0,0)-\int_{[0,x]}y_{3}(x-t, \lambda,0,0)y_{j}(t.\lambda,p,q)\mathrm{d}\tilde{\mu}_{p,q}(t)\\
\!\!\!\!&&\!\!\!\!-\int_{[0,x]}y_{3}(x-t,\lambda,0,0)2q(t)y_{j}^{\prime}(t,\lambda,p,q)\mathrm{d}t.
\end{eqnarray*}
Then using a variant of the integration by parts formula for the product of three functions, we have
\begin{eqnarray*}
y_{j}(x,\lambda,p,q)\!\!\!\!&=&\!\!\!\!y_{j}(x,\lambda, 0,0)-2q(t)y_{3}(x-t,\lambda,p,q)y_{j}(t,\lambda,p,q)\vert_{t=0}^{t=x+}\\
\!\!\!\!&&\!\!\!\!+\int_{[0,x]}y_{j}(t.\lambda,p,q)\mathrm{d}(y_{3}(x-t,\lambda,0,0)2q(t))\\
\!\!\!\!&&\!\!\!\!-\int_{[0,x]}y_{3}(x-t, \lambda,0,0)y_{j}(t,\lambda,p,q)\mathrm{d}\tilde{\mu}_{p,q}(t)\\
\!\!\!\!&=&\!\!\!\!y_{j}(x,\lambda, 0,0)-2q(t)y_{3}(x-t,\lambda,p,q)y_{j}(t,\lambda,p,q)\vert_{t=0}^{t=x+}\\
\!\!\!\!&&\!\!\!\!-\int_{[0,x]}y^{\prime}_{3}(x-t,\lambda,0,0)2q(t)y_{j}(t,\lambda,p,q)\mathrm{d}t\\
\!\!\!\!&&\!\!\!\!+\int_{[0,x]}y_{3}(x-t, \lambda,0,0)y_{j}(t,\lambda,p,q)\mathrm{d}{\mu}_{p,q}(t)\\
&=&\!\!\!\!y_{j}(x,\lambda ,0,0)-\int_{[0,x]}2q\left( t\right) y_{3}^{\prime
}\left( x-t,\lambda ,0,0\right) y_{j}\left( t,\lambda ,p,q\right) \mathrm{d}t\\
&&+\int_{[0,x]}y_{3}\left( x-t,\lambda ,0,0\right) y_{j}\left( t,\lambda
,p,q\right) \mathrm{d}\mu _{p,q}\left( t\right) ,
\end{eqnarray*}}
where $\mu _{p,q}\left( t\right) =q\left( t\right) +ip\left( t\right) $.
Following Picard's {iteration} we write%
\begin{equation}
y_{j}(x,\lambda ,p,q)=\underset{m\in {\mathbb{N}}_{0}}{\sum }c_{m}^{j}(x,\lambda
,p,q),  \label{ser}
\end{equation}%
where%
\begin{eqnarray*}
c_{0}^{j}(x,\lambda ,p,q)\!\!\!\! &=&\!\!\!\!y_{j}(x,\lambda ,0,0), \\
c_{m}^{j}(x,\lambda ,p,q)\!\!\!\! &=&\!\!\!\!-\int_{[0,x]}2q\left( t\right)
y_{3}^{\prime }\left( x-t,\lambda ,0,0\right) c_{m-1}^{j}({t},\lambda ,p,q)%
\mathrm{d}t \\
&&+\int_{[0,x]}y_{3}\left( x-t,\lambda ,0,0\right) c_{m-1}^{j}({t},\lambda
,p,q)\mathrm{d}\mu _{p,q}\left( t\right) \\
&=:&\!\!\!\!\int_{[0,x]}c_{m-1}^{j}({t},\lambda ,p,q)\left[ -2q\left(
t\right) y_{3}^{\prime }\left( x-t,\lambda ,0,0\right) \mathrm{d}t\right. \\
&&+\left. y_{3}\left( x-t,\lambda ,0,0\right) \mathrm{d}\mu _{p,q}\left(
t\right) \right] , m\in {\mathbb{N}}.
\end{eqnarray*}%
{Moreover, for $m\in \mathbb{N}$, it is easy to verify that}
\begin{eqnarray*}
c_{m}^{j}(x,\lambda ,p,q)\!\!\!\! \!\!&=&\!\!\!\!\!\!\int_{0\leqslant t_{1}<t_{2}<\cdots
t_{m}<t_{m+1}:=x}y_{j}(t_{1},\lambda ,0,0)\overset{m}{\underset{l=1}{\prod }}%
\left[ -2y_{3}^{\prime }\left( t_{l+1}-t_{l},\lambda
,0,0\right)\right. \\
&&\left. q\left( t_{l}\right)  \mathrm{d}t_{l}+y_{3}\left( t_{l+1}-t_{l},\lambda ,0,0\right) \mathrm{d}\mu
_{p,q}\left( t_{l}\right) \right] .
\end{eqnarray*}%
From (\ref{3}), we have
\begin{equation*}
\left\vert y_{3}\left( x,\lambda ,0,0\right) \right\vert \leqslant 3\Xi
\left( x,\lambda \right) , \left\vert y_{3}^{\prime }\left( x,\lambda
,0,0\right) \right\vert =\left\vert y_{2}\left( x,\lambda ,0,0\right)
\right\vert \leqslant 3\Xi \left( x,\lambda \right) .
\end{equation*}%
{Therefore,} in light of (\ref{i3}) and (\ref{222}), we have
\begin{eqnarray*}
\left\vert c_{m}^{j}(x,\lambda ,p,q)\right\vert \!\!\!\! &\leqslant&\!\!\!\!\int_{0\leqslant t_{1}<t_{2}<\cdots t_{m}<t_{m+1}:=x}\frac{3}{\left\vert k\right\vert ^{j-1}}\overset{m}{\underset{l=1}{\prod }}{3}\left[2\left\vert q\left( t_{l}\right) \right\vert \mathrm{d}t_{l}\right.\\
&&\left.+\mathrm{d}\mathbf{V}_{p}\left( t_{l}\right) +\mathrm{d}\mathbf{V}_{q}\left( t_{l}\right) \right]\Xi \left( t_{1},\lambda \right) \cdot \Xi \left( x-t_{1},\lambda \right)  \\
&=&\!\!\!\!\frac{1}{m!}\frac{3\Xi \left( x,\lambda \right) }{\left\vert
k\right\vert ^{j-1}}\left[ 3\int_{[0,x]}\left( 2\left\vert q\left( t\right)
\right\vert \mathrm{d}t+\mathrm{d}\mathbf{V}_{p}\left( t\right) +\mathrm{d}%
\mathbf{V}_{q}\left( t\right) \right) \right] ^{m} \\
&\leqslant &\!\!\!\!\frac{1}{m!}\frac{3\Xi \left( x,\lambda \right) }{\left\vert
k\right\vert ^{j-1}}\left[ 3\left( 2\Vert q \Vert_{\mathbf{V}}  +\mathbf{\check{p}}%
\left( x\right) +\mathbf{\check{q}}\left( x\right) \right) \right] ^{m}{,}
\end{eqnarray*}
{and thus}%
\begin{equation*}
\left\vert y_{j}(x,\lambda ,p,q)\right\vert \leqslant \frac{3}{\left\vert
k\right\vert ^{j-1}}\Xi \left( x,\lambda \right) e^{3\left( 2\Vert q \Vert_{\mathbf{V}}  +\mathbf{\check{p}}\left( x\right) +\mathbf{\check{q}}\left(
x\right) \right) }.
\end{equation*}%
Note that
\begin{eqnarray*}
\left\vert y_{3}^{\prime }\left( x,\lambda ,0,0\right) \right\vert \!\!\!\!
&=&\!\!\!\!\left\vert y_{2}\left( x,\lambda ,0,0\right) \right\vert
\leqslant \frac{3}{\left\vert k\right\vert }\Xi \left( x,\lambda \right) , \\
\left\vert y_{3}\left( x,\lambda ,0,0\right) \right\vert \!\!\!\!
&=&\!\!\!\!\left\vert \int_{[0,x]}y_{2}\left( t,\lambda ,0,0\right) \mathrm{d%
}t\right\vert \leqslant \frac{3}{\left\vert k\right\vert }\Xi \left(
x,\lambda \right) ,
\end{eqnarray*}%
then proceeding as {in} the proof of {the inequality} (\ref{ypq}), we obtain the inequality (\ref%
{pq-00}). \qed

\begin{rem}
\label{ori-e} {In fact, it is straightforward to show that the series in $(\ref{ser})$ converges uniformly for $x\in I$, $\lambda \in U$, $p \in B_{\delta_{p}}$ and $q \in B_{\delta_{q}}$, where $U$ is any bounded subset of ${\mathbb{C}}$, $B_{\delta_{p}}:=\{f \in {\mathcal{M}}_{0}(I,{\mathbb{R}}), \Vert f-p \Vert _{\mathbf{V}}\leqslant\delta_{p}\}$, $B_{\delta_{q}}:=\{f \in {\mathcal{M}}_{0}(I,{\mathbb{R}}), \Vert f-q \Vert _{\mathbf{V}}\leqslant\delta_{p}\}$, $\delta_{p} > 0$, $\delta_{q} > 0$. When $\frac{\mathrm{d}}{\mathrm{d}x}p\in {\mathcal{L}}^{2}(I,{\mathbb{R}})$, $q\in {\mathcal{H}}^{1}(I,{\mathbb{R}})$, Amour L gave the similar estimates for fundamental solutions of {the} equation $(\ref{equation})$ in $\cite[Theorem\ 2.4-2.5]{amour1999siam}$.}
\end{rem}

\begin{lem}
\label{entire} {For $(x,\lambda ,p,q)\in I \times {%
\mathbb{C\times }\mathcal{M}}_{0}(I,{\mathbb{R}})\times {\mathcal{M}}_{0}(I,{%
\mathbb{R}})$, $j=1, 2, 3$, we have%
\begin{equation*}
y_{j}(x,\lambda ,p,q),\quad y_{j}^{\prime }(x,\lambda ,p,q),\quad \left(
y_{j}^{\prime }\right) ^{\bullet }(x,\lambda ,p,q)
\end{equation*}%
are entire functions of $\lambda $.}
\end{lem}

\Proof The proof is similar to that of \cite[Theorem
2.6]{amour1999siam}. \qed

\section{Eigenvalue of Measure Differential Equation}
This section is devoted to study the eigenvalues of the boundary value problems (\ref{equation})-$(BC)_{\xi }$, $\xi =1,2$ with coefficients $p$, $q\in {\mathcal{M}}_{0}(I,{\mathbb{R}})$.
\subsection{The Distribution of Eigenvalues}
In this subsection, we investigate the counting lemma {(see Theorem \ref{countinglemma})} for the boundary value problems (\ref{equation})-$(BC)_{\xi }$, $\xi =1,2$, which implies the distribution and estimates of eigenvalues. Firstly, we give some notations and basic lemmas.
\begin{defn}
\label{ev-def} {For $p$, $q \in {\mathcal{M}}_{0}(I,{\mathbb{R}})$, a complex number $\lambda $ is called an eigenvalue of the boundary value problem $(\ref{equation})$-$\left( BC\right)_{1} $ if the equation $(\ref{equation})$ with such a parameter $\lambda $ has a nontrivial solution $e(x,\lambda ,p,q)$ on $I$ satisfying the boundary conditions $\left( BC\right)_{1} $. The solution $e(x,\lambda ,p,q) $ is called an eigenfunction of $\lambda $. The number of {linearly} independent eigenfunctions associated with $\lambda $ is called the geometric multiplicity $($g$-multiplicity)$ of $\lambda $. The eigenvalues and eigenfunctions for the boundary value problem $(\ref{equation})$-$(BC)_{2}$ are defined similarly. }
\end{defn}

\begin{lem}
\label{ev-real}{For $p$, $q\in {\mathcal{M}}_{0}(I,{\mathbb{R}})$, all eigenvalues of {the} boundary value problems $(\ref{equation})$-$(BC)_{\xi }$, $\xi=1,2$ are real.}
\end{lem}

\Proof Suppose $\lambda$, $\mathrm{Im} \lambda \neq 0$, is an eigenvalue of the boundary value problem (\ref{equation})-$\left( BC\right)_{1} $, then the corresponding eigenfunction {$e:=e(x,\lambda, p,q)$} satisfies
\begin{equation}
i\mathrm{d}\left( e^{\prime }\right) ^{\bullet }+2iq\left(
x\right) e^{\prime }\mathrm{d}x+e\left( i\mathrm{d}q\left( x\right) +\mathrm{d}p\left(
x\right) \right) = \lambda e\mathrm{d}x,  \label{yuan}
\end{equation}
and
\begin{equation}
-i\mathrm{d}\left( \bar{e}^{\prime }\right) ^{\bullet }-2iq\left(
x\right) \bar{e}^{\prime }\mathrm{d}x+\bar{e}\left(- i\mathrm{d}q\left( x\right) +\mathrm{d}p\left(
x\right) \right) = \bar{\lambda} \bar{e}\mathrm{d}x.  \label{baryuan}
\end{equation}
Here, $\bar{\alpha}$ denotes the {conjugation} of $\alpha$. {Multiplying (\ref{yuan}) by $\bar{e}$, (\ref{baryuan}) by ${e}$, and taking the difference, we find
\begin{equation*}
i\mathrm{d}\left( e^{\prime }\right) ^{\bullet }\bar{e} +i\mathrm{d}\left( \bar{e}^{\prime }\right) ^{\bullet }e+2iq\left(
x\right) e^{\prime }\bar{e}\mathrm{d}x+2iq\left(x\right) \bar{e}^{\prime }e\mathrm{d}x+2ie\bar{e}\mathrm{d}q\left( x\right)= 2i\mathrm{Im}{\lambda} e\bar{e}\mathrm{d}x.
\end{equation*}
{Hence,
\begin{eqnarray*}
\int_{[0,1]}i\bar{e}\mathrm{d}\left( e^{\prime }\right) ^{\bullet } +\int_{[0,1]}ie\mathrm{d}\left( \bar{e}^{\prime }\right) ^{\bullet }+\int_{[0,1]}2iq\left(x\right) e^{\prime }\bar{e}\mathrm{d}x+\int_{[0,1]}2iq\left(x\right) \bar{e}^{\prime }e\mathrm{d}x &&\\
+\int_{[0,1]}2ie\bar{e}\mathrm{d}q\left( x\right)= \int_{[0,1]}2i\mathrm{Im}{\lambda} e\bar{e}\mathrm{d}x.&&
\end{eqnarray*}}
Using the integration by parts formula, we have}
\begin{equation*}
\left.(i\bar{e} (e^{\prime })^{\bullet }+i(\bar{e}^{\prime })^{\bullet } e-i\bar{e}^{\prime } e^{\prime}+2iq(x)e\bar{e}) \right\vert _{x=0}^{x=1+}=2i\mathrm{Im} \lambda\int_{I} \vert e\vert ^{2}\mathrm{d}x.
\end{equation*}
According to the boundary conditions $\left( BC\right)_{1} $, $q(0)=0$, and $\int_{I} \vert e\vert ^{2}\mathrm{d}x \neq 0$, we obtain that $\mathrm{Im} \lambda = 0$. Similarly, we can prove that the eigenvalues of the boundary value problem (\ref{equation})-$(BC)_{2}$ are all real. \qed

\begin{lem}
\label{doe}{Fix $(p,q)\in {\mathcal{M}}_{0}(I,{\mathbb{R}})\times {%
\mathcal{M}}_{0}(I,{\mathbb{R}})$, and let
\begin{equation*}
M_{\xi }\left( \lambda ,p,q\right) =\left(
\begin{array}{cc}
y_{1}(1,\lambda ,p,q) & y_{2}(1,\lambda ,p,q) \\
y_{1}^{\prime }(1,\lambda ,p,q) & y_{2}^{\prime }(1,\lambda ,p,q)+(-1)^{\xi}%
\end{array}%
\right) , \quad \xi =1,2.
\end{equation*}
}

{$(i)$ For $\xi =1,2$, {each} eigenvalue $\lambda$ of {the} boundary value problem $(\ref{equation})$-$\left( BC\right)_{\xi} $ is of $g$-multiplicity one or two and it is {a} root of%
\begin{equation*}
\Delta _{\xi }\left( \lambda \right) :=\det M_{\xi }\left( \lambda
,p,q\right).
\end{equation*}
}

{$(ii)$ For $\xi =1,2$, suppose $\lambda $ is an eigenvalue of
{the} boundary value problem $(\ref{equation})$-$(BC)_{\xi}$, then the $g$-multiplicity of $\lambda $ is two if and only if%
\begin{equation*}
M_{\xi }\left( \lambda ,p,q\right) =\left(
\begin{array}{cc}
0 & 0 \\
0 & 0%
\end{array}%
\right) .
\end{equation*}}

{$(iii)$ Recall the notations in $(\ref{L1})$. If $\Lambda _{1}\left( p,q\right) \cap \Lambda _{2}\left(
p,q\right) =\emptyset $, then the eigenvalues of {the} boundary value problems $(\ref{equation})$-$(BC)_{\xi }$, $\xi =1,2$ are of $g$-multiplicity one.}
\end{lem}

\Proof The proof{s} of (i) and (ii) {are} similar to {those} of \cite[Theorem 3.1]{amour1999siam}.

(iii) From (i) {and} a simple calculation{, we infer that}
\begin{eqnarray*}
\Delta _{\xi }\left( \lambda \right)\hh &=&\hh y_{1}(1,\lambda,p,q)y_{2}^{\prime }(1,\lambda,p,q)-y_{1}^{\prime }(1,\lambda,p,q)y_{2}(1,\lambda,p,q)+(-1)^{\xi} y_{1}(1,\lambda,p,q)\\
&=&\hh [y_{2}(1,\lambda,p,q),y_{1}(1,\lambda,p,q)]+(-1)^{\xi}y_{1}(1,\lambda,p,q),\quad  \xi =1,2.
\end{eqnarray*}
{By the equality}
\begin{equation*}
[y_{2}(1,\lambda,p,q),y_{1}(1,\lambda,p,q)]=\bar{y}_{1}(1,\bar{\lambda},p,q)
\end{equation*}
given by Mckean H P \cite[p. 614]{mckean1981cpam}, we have
\begin{equation}
\Delta _{\xi }\left( \lambda \right) =\bar{y}_{1}(1,\bar{\lambda},p,q)+(-1)^{\xi}
y_{1}(1,\lambda ,p,q),  \quad \xi =1,2.\label{D}
\end{equation}
Apparently, $\bar{y}_{1}({x},\bar{\lambda},p,q)$ is the solution of {the} equation
\begin{equation*}
i\mathrm{d}\left( y^{\prime }\right) ^{\bullet }+2iq\left(
x\right) y^{\prime }\mathrm{d}x+y\left( i\mathrm{d}q\left( x\right) -\mathrm{d}p\left(
x\right) \right) = -\lambda y\mathrm{d}x, \qquad x\in (0,1)
\end{equation*}
with initial conditions $(\bar{y}_{1}(0,\bar{\lambda},p,q),\bar{y}_{1}^{\prime }(0,\bar{\lambda},p,q),(\bar{y}_{1}^{\prime })^{\bullet }(0,\bar{\lambda},p,q))=(1,0,0)$.
For each $\lambda \in {\mathbb{R}}$, the following identities
\begin{eqnarray*}
\Delta _{1}\left( \lambda \right) \!\!\!\! &=&\!\!\!\!\bar{y}_{1}(1,\lambda
,p,q)-y_{1}(1,\lambda ,p,q)=-2i\mathrm{Im}y_{1}(1,\lambda ,p,q), \\
\Delta _{2}\left( \lambda \right) \!\!\!\! &=&\!\!\!\!\bar{y}_{1}(1,\lambda
,p,q)+y_{1}(1,\lambda ,p,q)=2\mathrm{Re}y_{1}(1,\lambda ,p,q)
\end{eqnarray*}%
hold. Since the eigenvalues of {the} boundary value problems (\ref{equation})-$(BC)_{\xi }$, $\xi =1,2$ are real and $\Lambda
_{1}\left( p,q\right) \cap \Lambda _{2}\left( p,q\right) =\emptyset $, we
obtain that {as a function of $\lambda$, }$y_{1}(1,\lambda ,p,q)$ has no zeros in ${\mathbb{R}}${. Hence,}%
\begin{equation*}
M_{\xi }\left( \lambda ,p,q\right) \neq \left(
\begin{array}{cc}
0 & 0 \\
0 & 0%
\end{array}%
\right) ,  \xi =1,2,  \lambda \in {\mathbb{R}}.
\end{equation*}%
{Then the statement $(iii)$ of this lemma follows from the statement $(ii)$ of this lemma.}\qed

\begin{defn}
{For $\xi =1,2$, the order of an eigenvalue $\lambda$ as a
root of $\Delta _{\xi }\left( \lambda \right) =0$ is called the
algebraic multiplicity $($$a$-multiplicity$)$ of $\lambda$.}
\end{defn}

\begin{lem}
\label{omiga}
{Denote $\Omega _{\xi }:=\overset{3}{\underset{j=1}{\cap }}\Omega
_{\xi ,j}$, $\xi =1,2$, where
\begin{equation*}
\Omega _{\xi ,j}:= \left\{ k\in {\mathbb{C}}; \left\vert
\omega ^{j-1}k-(2n+\xi -1)\pi \right\vert \geqslant \frac{\pi }{6}, n\in \mathbb{Z}\right\} ,j=1,2,3, {\omega=e^{\frac{2\pi i}{3}}}.
\end{equation*}
Then for $k\in \Omega _{\xi }$, $\xi =1,2$,
there is a constant $C_{\pi }>0$, which is independent of $j$, $\omega$ and $k${,} such that
\begin{equation*}
e^{\left\vert \mathrm{Im}\frac{\omega ^{j-1}k}{2}\right\vert }<C_{\pi }^{%
\frac{1}{3}}\left\vert \sin \frac{\omega ^{j-1}k}{2}\right\vert ,\
e^{\left\vert \mathrm{Im}\frac{\omega ^{j-1}k}{2}\right\vert }<C_{\pi }^{%
\frac{1}{3}}\left\vert \cos \frac{\omega ^{j-1}k}{2}\right\vert , j=1,2,3.
\end{equation*}%
}
\end{lem}
\Proof See \cite[p. 27]{ptrbookisp}. \qed

Now we give the main result of this subsection{. We mention that the following result gives} an explanation of the indexation $n$ of the eigenvalue $\lambda_{\xi, n} (p,q)$, $\xi =1,2$.

\begin{thm}
\label{countinglemma}{$\left( the\text{ }counting\text{ }lemma\right) $ Suppose $(p,q)\in {\mathcal{M}}_{0}(I,{\mathbb{R}})\times {%
\mathcal{M}}_{0}(I,{\mathbb{R}})$.}

{$\left( i\right) $ Let $N$ be an integer satisfying
\begin{equation*}
\left( 2N+1\right) \pi >\frac{9C_{\pi }}{4}e^{3\left( 3\Vert q \Vert_{\mathbf{V}}+\Vert p \Vert_{\mathbf{V}} \right)}.
\end{equation*}
Then the boundary value problem $(\ref{equation})$-$(BC)_{1}$ has exactly %
$2N+1$ eigenvalues, counted with $a$-multiplicities, in the open $\lambda $%
-disc%
\begin{equation*}
\left\{ \lambda =k^{3}\in {\mathbb{C}}; \left\vert k\right\vert <\left(
2N+1\right) \pi \right\}{,}
\end{equation*}%
and exactly one {algebraically simple} eigenvalue in each open $\lambda $-disc%
\begin{equation*}
\left\{ \lambda =k^{3}\in {\mathbb{C}}; \left\vert k-2n\pi \right\vert <%
\frac{\pi }{3}\right\}
\end{equation*}%
for $\left\vert n\right\vert \geqslant N$.}

{$\left( ii\right) $ Let $N$ be an integer satisfying
\begin{equation*}
2N\pi >\max \left\{ \frac{9C_{\pi }}{2}e^{3\left( 3\Vert q \Vert_{\mathbf{V}}+\Vert p \Vert_{\mathbf{V}} \right) },2\ln\frac{C_{\pi }}{2}, 2\sqrt{2}\ln \frac{C_{\pi }}{4}\right\} .
\end{equation*}
Then the boundary value problem $(\ref{equation})$-$(BC)_{2}$ has exactly $2N$ eigenvalues, counted with
$a$-multiplicities, in the open $\lambda $-disc%
\begin{equation*}
\left\{ \lambda =k^{3}\in {\mathbb{C}}; \left\vert k\right\vert <2N\pi
\right\}{,}
\end{equation*}%
and exactly one {algebraically simple} eigenvalue in each open $\lambda $-disc%
\begin{equation*}
\left\{ \lambda =k^{3}\in {\mathbb{C}}; \left\vert k-\left( 2n+1\right) \pi
\right\vert <\frac{\pi }{3}\right\}
\end{equation*}%
for $\left\vert n\right\vert \geqslant N$.}
\end{thm}

\Proof We divide our proof {into} two steps.

\textit{Step1}. For each $\lambda \in {\mathbb{C}}$, $p$, $q\in {\mathcal{M}}_{0}(I,{\mathbb{R}})$, let $ \left(
\begin{array}{l}
Y_{1}(x,\lambda ,p,q) \\
Z_{1}(x,\lambda ,p,q)%
\end{array}%
\right) $ denote the solution of {the} equation
\begin{equation*}
\left(
\begin{array}{l}
-\mathrm{d}\left( Z^{\prime }\right) ^{\bullet}-2q\left( x\right) Z^{\prime }\mathrm{d}x-Z\mathrm{d}q\left( x\right)+Y\mathrm{d}p\left( x\right)  \\
\mathrm{d}\left( Y^{\prime }\right) ^{\bullet }+2q\left( x\right) Y^{\prime }\mathrm{d}x+Y\mathrm{%
d}q\left( x\right) +Z\mathrm{d}p\left( x\right)%
\end{array}%
\right)=\lambda \left(
\begin{array}{l}
Y\mathrm{d}x \\
Z\mathrm{d}x%
\end{array}%
\right)
 , x\in I,
\end{equation*}%
with {the} initial conditions%
\begin{equation*}
\left(
\begin{array}{l}
Y_{1}(0,\lambda ,p,q) \\
Z_{1}(0,\lambda ,p,q)%
\end{array}%
\right) =\left(
\begin{array}{l}
1 \\
0%
\end{array}%
\right) , \left(
\begin{array}{l}
Y_{1}^{\prime }(0,\lambda ,p,q) \\
Z_{1}^{\prime }(0,\lambda ,p,q)%
\end{array}%
\right) =\left(
\begin{array}{l}
0 \\
0%
\end{array}%
\right) , \left(
\begin{array}{l}
\left( Y_{1}^{\prime }\right) ^{\bullet }(0,\lambda ,p,q) \\
\left( Z_{1}^{\prime }\right) ^{\bullet }(0,\lambda ,p,q)%
\end{array}%
\right) =\left(
\begin{array}{l}
0 \\
0%
\end{array}%
\right) .
\end{equation*}%
{Here,} $Y_{1}(x,\lambda ,p,q)$ and $Z_{1}(x,\lambda ,p,q)$ are real-valued for $\lambda \in {\mathbb{R}}$. For $x\in I${,} $\lambda \in {%
\mathbb{C}}$, {and $k=\lambda^{\frac{1}{3}}$, }a straightforward calculation gives
\begin{eqnarray}
y_{1}(x,\lambda ,p,q)\hh&=&\hh Y_{1}(x,\lambda ,p,q)+iZ_{1}(x,\lambda ,p,q), \label{Y+Z}\\
\bar{y}_{1}(x,\bar{\lambda},p,q)\hh&=&\hh Y_{1}(x,\lambda ,p,q)-iZ_{1}(x,\lambda
,p,q),  \label{Y-Z}\\
Y_{1}(x,\lambda ,0,0)\hh&=&\hh \frac{1}{3}\left( \cos k{x}+\cos \omega k{x}+\cos \omega
^{2}k{x}\right)\notag\\
&=&\hh\frac{1}{3}%
\left( 4\cos \left( \frac{k{x}}{2}\right) \cos \left( \frac{\omega k{x}}{2}%
\right) \cos \left( \frac{\omega ^{2}k{x}}{2}\right) -1\right) , \notag\\
Z_{1}(x,\lambda ,0,0)\hh&=&\hh \frac{1}{3}\left( \sin k{x}+\sin \omega k{x}+\sin \omega
^{2}k{x}\right)\notag \\
&=&\hh{-}\frac{4}{3%
}\sin \left( \frac{k{x}}{2}\right) \sin \left( \frac{\omega k{x}}{2}\right) \sin
\left( \frac{\omega ^{2}k{x}}{2}\right) .\notag
\end{eqnarray}
According to (\ref{D})-(\ref{Y-Z}), it follows that for $\lambda \in \mathbb{C}$,
\begin{equation*}
\Delta _{1}\left( \lambda \right) =-2iZ_{1}(1,\lambda ,p,q), \quad \Delta _{2}\left( \lambda \right) =2Y_{1}(1,\lambda ,p,q).
\end{equation*}%
Then {by} Lemma \ref{doe} (i), {we know that} in order to prove Theorem \ref{countinglemma}, it is sufficient to discuss the zeros of $Z_{1}(1,\lambda ,p,q)$ and $Y_{1}(1,\lambda ,p,q)$, respectively.

\textit{Step 2}. {In view of} Theorem \ref{eopq} and Lemma \ref{entire}, {an argument similar to the one used} in \cite[Appendix]{amgu} {shows} that $Y_{1}(1,\lambda ,p,q)$, $Z_{1}(1,\lambda ,p,q)$ are entire functions of $\lambda \in \mathbb{C}${,} and
\begin{eqnarray}
&&\left\vert Z_{1}(x,\lambda ,p,q)-Z_{1}(x,\lambda ,0,0)\right\vert
\leqslant \frac{3}{\left\vert k\right\vert }\Xi \left( x,\lambda \right)
e^{3\left( 2\Vert q \Vert_{\mathbf{V}} +\mathbf{\check{p}}\left( x\right) +\mathbf{%
\check{q}}\left( x\right) \right) } , \label{D+-D0} \\
&&\left\vert Y_{1}(x,\lambda ,p,q)-Y_{1}(x,\lambda ,0,0)\right\vert
\leqslant \frac{3}{\left\vert k\right\vert }\Xi \left( x,\lambda \right)
e^{3\left( 2\Vert q \Vert_{\mathbf{V}} +\mathbf{\check{p}}\left( x\right) +\mathbf{%
\check{q}}\left( x\right) \right) }\label{D--D0}
\end{eqnarray}
hold for $(x,\lambda ,p,q)\in I \times {\mathbb{C}\times \mathcal{M}}_{0}(I,{\mathbb{R}})\times {\mathcal{M}}_{0}(I,{\mathbb{R}})$. {In view of} Lemma \ref{omiga} {and} \cite[Lemma 3.5]{amour1999siam}, we get%
\begin{eqnarray}
&&\left\vert Z_{1}(1,\lambda ,0,0)\right\vert >\frac{4}{3C_{\pi }}%
\Xi \left( 1,\lambda \right) , \forall k\in \Omega _{1},  \label{D0+} \\
&&\left\vert Y_{1}(1,\lambda ,0,0)\right\vert >\frac{2}{3C_{\pi }}%
\Xi \left( 1,\lambda \right) , \forall k\in \Omega _{2}\backslash \left\{
k\in {\mathbb{C}}; \left\vert k\right\vert <2\ln \frac{C_{\pi }}{2}\right\}
.  \label{D0-}
\end{eqnarray}
According to the inequalities (\ref{D+-D0}) and (\ref{D0+}%
), it follows that
\begin{equation*}
\left\vert Z_{1}(1,\lambda ,p,q) -Z_{1}(1,\lambda ,0,0) \right\vert < \vert Z_{1}(1,\lambda ,0,0) \vert
\end{equation*}
holds for $k\in \Omega _{1}$ satisfying $\vert k \vert >\frac{9C_{\pi }}{4}e^{3\left( 3\Vert q \Vert_{\mathbf{V}}+\Vert p \Vert_{\mathbf{V}} \right)
}$. {Now} we select an integer $N$ satisfying $\left( 2N+1\right) \pi >\frac{9C_{\pi }}{4}e^{3\left( 3\Vert q \Vert_{\mathbf{V}}+\Vert p \Vert_{\mathbf{V}} \right)}${. Then} using Rouch\'{e} theorem{,} we obtain {that} $Z_{1}(1,\lambda ,p,q)$ and $ Z_{1}(1,\lambda ,0,0)$ have the same number of zeros in the $\lambda $-discs defined in (i). Since {$ Z_{1}(1,\lambda ,0,0)$ has only the simple zeros} $\lambda _{1,n} ^{0} =(2n\pi)^{3}$, {$n\in \mathbb{Z}$, the statement in $(i)$ follows}. %

It remains to characterize the distribution of zeros of $Y_{1}(1,\lambda,p,q)$. Let
\begin{equation*}
Y_{1}^{0}(\lambda):=Y_{1}(1,\lambda,0,0)+\frac{1}{3}=\frac{4}{3}\cos \left( \frac{k}{2}\right) \cos \left( \frac{\omega k}{2}\right) \cos \left( \frac{\omega ^{2}k}{2}\right).
\end{equation*}
{Then} from Lemma \ref{omiga}, we have
\begin{equation*}
\vert Y_{1}^{0}(\lambda)\vert>\frac{4}{3C_{\pi }}\Xi \left( 1,\lambda \right) \geqslant \frac{4}{3C_{\pi }}e^{\frac{\sqrt{2}}{4}\vert k\vert} , \forall k\in \Omega _{2}{.}
\end{equation*}
{Therefore,} for any $k\in \Omega _{2}$ satisfying $\vert k\vert>2\sqrt{2}\ln \frac{C_{\pi }}{4}$, {it is easy to see that}
\begin{equation*}
\vert Y_{1}(1,\lambda,0,0)-Y_{1}^{0}(\lambda)\vert<\vert Y_{1}^{0}(\lambda)\vert.
\end{equation*}
{Now} let $N$ be an integer satisfying $2N\pi >2\sqrt{2}\ln \frac{C_{\pi }}{4}$. {Then} using Rouch\'{e} theorem, we obtain that $Y_{1}(1,\lambda,0,0) $ and $Y_{1}^{0}(\lambda)$ have the same number of zeros in {the} $\lambda$-discs
\begin{equation}
\left\{ \lambda =k^{3}\in {\mathbb{C}}; \left\vert k\right\vert <2N\pi
\right\} \label{2N}
\end{equation}
and
\begin{equation}
\left\{ \lambda =k^{3}\in {\mathbb{C}}; \left\vert k-\left( 2n+1\right) \pi
\right\vert <\frac{\pi }{3}\right\}, \vert n\vert \geqslant N. \label{NN}
\end{equation}
{Combining} (\ref{D--D0}) and (\ref{D0-}), we obtain the inequality
\begin{equation*}
\left\vert Y_{1}(1,\lambda,p,q) -Y_{1}(1,\lambda,0,0) \right\vert < \vert Y_{1}(1,\lambda,0,0)\vert
\end{equation*}
is true for $k \in \Omega _{2}$ satisfying $\vert k \vert >\max \left\{ \frac{9C_{\pi }}{2}e^{3\left( 3\Vert q \Vert_{\mathbf{V}}+\Vert p \Vert_{\mathbf{V}} \right) },2\ln
\frac{C_{\pi }}{2}\right\} $. {Thus} for any integer $N$ satisfying $2N\pi >\max \left\{ 2\ln\frac{C_{\pi }}{2}, \frac{9C_{\pi }}{2}e^{3\left( 3\Vert q \Vert_{\mathbf{V}}+\Vert p \Vert_{\mathbf{V}} \right) },2\sqrt{2}\ln \frac{C_{\pi }}{4}\right\} $, using Rouch\'{e} theorem again, we obtain that $Y_{1}(1,\lambda,0,0) $ and $Y_{1}(1,\lambda,p,q)$ have the same number of zeros in {the} $\lambda$-discs (\ref{2N}) and (\ref{NN}). Hence{,} $Y_{1}^{0}(\lambda)$ and $Y_{1}(1,\lambda,p,q)$ have the same number of zeros in {the} $\lambda$-discs (\ref{2N}) and (\ref{NN}). Since the zeros of {the} entire function $Y_{1}^{0}(\lambda)$ are $\lambda _{2,n} ^{0} =((2n+1)\pi)^{3}$, {$n\in\mathbb{Z}$}, we obtain the statement in (ii). \qed

As a consequence of Theorem \ref{countinglemma}, the following {result} gives a rough asymptotic expansion of the eigenvalues of the boundary value problems (\ref{equation})-$(BC)_{\xi }$, $\xi =1,2$.

\begin{cor}
\label{estimate}
{For $(p,q)\in {\mathcal{M}}_{0}(I,{\mathbb{R}})\times {\mathcal{M}}%
_{0}(I,{\mathbb{R}})$, we have
\begin{equation*}
\lambda _{1,n}\left( p,q\right)=(2n\pi )^{3}-4n\pi \int_I  q(x) \mathrm{d}x+O(1),
\end{equation*}
and
\begin{equation*}
\lambda _{2,n}\left( p,q\right)=((2n+1)\pi )^{3}-2(2n+1)\pi \int_I  q(x) \mathrm{d}x+O(1)
\end{equation*}
as $\left \vert n \right \vert \rightarrow +\infty $.}
\end{cor}

\Proof The proof is similar to {that of} \cite[Theorem 1.2]{amour1999siam}. \qed

\begin{cor}
\label{simple} {Let
\begin{eqnarray*}
N_{0}\hh&:=&\hh\min \left\{ N\in \mathbb{N} \bigg |
\left( 2N+1\right) \pi >\frac{9C_{\pi }}{4}e^{3\left( 3\Vert q \Vert_{\mathbf{V}}+\Vert p \Vert_{\mathbf{V}} \right)}, \right.\\
&&\hh\left.2N\pi >\max \left\{ \frac{9C_{\pi }}{2}e^{3\left( 3\Vert q \Vert_{\mathbf{V}}+\Vert p \Vert_{\mathbf{V}} \right) },2\ln\frac{C_{\pi }}{2}, 2\sqrt{2}\ln \frac{C_{\pi }}{4}\right\}\right\} ,
\end{eqnarray*}%
then the $g$-multiplicity and $a$-multiplicity of each {eigenvalue} $\lambda _{\xi ,n}\left(
p,q\right) $, $\xi =1,2$, $\left\vert n\right\vert \geqslant N_{0}$, are equal to one. }
\end{cor}
\Proof {In view of} Theorem $\ref{countinglemma}${,} we obtain %
\begin{equation*}
\left\{ \lambda _{1 ,n}\left( p,q\right) ;n\geqslant N_{0}\right\} \cap \left\{ \lambda _{2,n}\left( p,q\right)
;n\geqslant N_{0}\right\} =\emptyset{,}
\end{equation*}%
and the $a$-multiplicity of each {eigenvalue} $\lambda _{\xi ,n}\left(
p,q\right) $, $\xi =1,2$, $\left\vert n\right\vert \geqslant N_{0}${,} is one. From Lemma $\mathrm{\ref{doe}}$ $(iii)$, we deduce that for any $\left\vert n\right\vert \geqslant N_{0}$, $\xi =1,2$, the
$g$-multiplicity of $\lambda _{\xi ,n}\left(
p,q\right) $ also equals one. \qed

\subsection{Dependence of Eigenvalues on Measures $p$, $q$}

In this subsection, we give the proofs of Theorem \ref{scoe} and Theorem \ref{frechet} announced in the introduction.

\smallskip\smallskip\noindent\emph{Proof of Theorem \ref{scoe}.} (i) Firstly, we discuss the dependence of {the} eigenvalues $\lambda_{\xi,n}(p,q)$, $\xi=1,2$ on {the measure} $q\in {({\mathcal{M}}_{0}(I,{\mathbb{R}}),w^{\ast })}$. Suppose that the sequence $\{q_{m}\}_{m \in \mathbb{N}}$ converges to $ q_{0}$ in ${({\mathcal{M}}_{0}(I,{\mathbb{R}}),w^{\ast })}$, then {from Lemma \ref{ubp},} there exists a constant $C_{q}>0$ such that $\underset{m\in {\mathbb{N}_{0}}}\sup\Vert q_{m}\Vert _{\mathbf{V}}\leqslant C_{q}$. Let
\begin{eqnarray*}
N_{1}\hh&:=&\hh\min \left\{ N\in {\mathbb{N}}\bigg |\left( 2N+1\right) \pi >\frac{9C_{\pi }}{4}e^{3\left(
3 C_{q}+\Vert p \Vert_{\mathbf{V}} \right) },\right.\\
&&\left.2N\pi >\max
\left\{ \frac{9C_{\pi }}{2}e^{3\left( 3 C_{q}+\Vert p \Vert_{\mathbf{V}} \right) },2\ln \frac{C_{\pi }}{2}, 2\sqrt{2}\ln \frac{C_{\pi }}{4}\right\} \right\}.
\end{eqnarray*}
For any $m\in \mathbb{N}_{0}$, {it follows from} Lemma \ref{entire} {that} $\Delta_{1}(\lambda ,q_{m}):={\bar{y}_{1}\left( 1,\bar{\lambda },p,q_{m}\right)-y_{1}\left( 1,\lambda ,p,q_{m}\right)}$ is an entire function of $\lambda\in \mathbb{C}$. For any integer $N \geqslant N_{1}$, let
\begin{equation*}
\gamma_{n}:=\{\lambda:\vert\lambda-\lambda_{1,n}(p,q_{0})\vert=\epsilon\} \mbox{ and }
\Gamma_{N}:=\overset{n=N}{\underset{n=-N}{\cup}}\gamma_{n},
\end{equation*}
where $\epsilon >0$ is any sufficiently small constant such that the contours $\gamma _{n}$, $\vert n\vert \leqslant N$, are disjoint and $\Delta_{1}(\lambda ,q_{0})\neq 0$ on $\Gamma_{N}$. Hence, there exists a constant $C_{q_{0}}>0$ such that
\begin{equation}
\vert \Delta_{1}(\lambda ,q_{0})\vert >C_{q_{0}}>0 \mbox{ on } \Gamma_{{N}}. \label{cq0}
\end{equation}
On the other hand, {in view of} {Remark \ref{uc1}}, we deduce that {as $m$ tends to infinity,} the sequence $\{\Delta_{1}(\lambda ,q_{m})\}$ converges to $\Delta_{1}(\lambda ,q_{0})$ uniformly on $\Gamma_{N}$. This implies that there exists a constant $M_{N}>0$ such that if $m>M_{N}$, one has
\begin{equation}
\vert\Delta_{1}(\lambda ,q_{m})-\Delta_{1}(\lambda ,q_{0})\vert<C_{q_{0}} \mbox{ on } \Gamma_{N}. \label{cqm0}
\end{equation}
{Therefore, combining} (\ref{cq0}) and (\ref{cqm0}), one deduces that for $m>M_{N}$,
\begin{equation*}
\vert\Delta_{1}(\lambda ,q_{m})-\Delta_{1}(\lambda ,q_{0})\vert<\vert \Delta_{1}(\lambda ,q_{0})\vert  \mbox{ on } \Gamma_{N}.
\end{equation*}
Then {by} Rouch\'{e} theorem, we {see that} for $m>M_{N}$, $\Delta_{1}(\lambda ,q_{m})$ and $\Delta_{1}(\lambda ,q_{0})$ have the same number of zeros inside each contour $\gamma_{n}$. Additionally, {in view of} Theorem \ref{countinglemma}, {we obtain that in the $\lambda$-disc $\left\{ \lambda =k^{3}\in {\mathbb{C}}; \left\vert k\right\vert <\left( 2N+1\right) \pi \right\} $, $\Delta_{1}(\lambda ,q_{0})$ has exactly $2N+1$ zeros $\lambda_{1,n}(p,q_{0})$, $\vert n\vert \leqslant N$, and $\Delta_{1}(\lambda ,q_{m})$ has exactly $2N+1$ zeros $\lambda_{1,n}(p,q_{m})$, $\vert n\vert \leqslant N$}. Hence, we obtain that given any {sufficiently} small $\epsilon >0$, there exists a constant $M_{N}>0$ such that for $m>M_{N}$,
\begin{equation*}
\vert \lambda_{1,n}(p,q_{m})-\lambda_{1,n}(p,q_{0}) \vert <\epsilon, \vert n\vert\leqslant N.
\end{equation*}
{Therefore,} from the arbitrariness of $N$, we obtain that each eigenvalue $\lambda_{1,n}(p,q)$ is continuous in $q \in {({\mathcal{M}}_{0}(I,{\mathbb{R}}),w^{\ast})}$.

Analogously, we can obtain $\lambda _{2,n}\left( p,q\right) $ is continuous in $q \in {({\mathcal{M}}_{0}(I,{\mathbb{R}}),w^{\ast})}$.

(ii) Now we deduce the dependence of {the} eigenvalues $\lambda_{\xi,n}(p,q)$, $\xi=1,2$, on {the measure} $p \in ({\mathcal{M}}_{0}(I,{\mathbb{R}}),w^{\ast})$. Suppose that the sequence $\{p_{m}\}_{m \in \mathbb{N}}$ converges to $p_{0}$ in $({\mathcal{M}}_{0}(I,{\mathbb{R}}),w^{\ast})$, then from Lemma \ref{ubp}, there exists a constant $C_{p_{0}}^{\ast}>0$ such that $\underset{m\in {\mathbb{N}_{0}}}\sup\Vert p_{m}\Vert _{\mathbf{V}}\leqslant C_{p_{0}}^{\ast}$. Let
\begin{eqnarray*}
N_{1}\hh&:=&\hh\min \left\{ N\in {\mathbb{N}}\bigg |\left( 2N+1\right) \pi >\frac{9C_{\pi }}{4}e^{3\left(
3\Vert q \Vert_{\mathbf{V}}+C_{p_{0}}^{\ast} \right) }, \right.\\
&&\left.2N\pi >\max\left\{ \frac{9C_{\pi }}{2}e^{3\left( 3\Vert q \Vert_{\mathbf{V}}+C_{p_{0}}^{\ast} \right) },2\ln \frac{%
C_{\pi }}{2}, 2\sqrt{2}\ln \frac{C_{\pi }}{4}\right\} \right\}.
\end{eqnarray*}
{In view of} {Remark \ref{uc1}}, we deduce that {as $m$ tends to infinity,} the sequence $\{\Delta_{1}(\lambda ,p_{m})\}$ converges to $\Delta_{1}(\lambda ,p_{0})$ uniformly on any bounded subset $U \subset \mathbb{C}$, where
\begin{equation*}
\Delta_{1}(\lambda ,p_{m})={\bar{y}_{1}\left( 1,\bar{\lambda },p_{m},q\right)-y_{1}\left( 1,\lambda ,p_{m},q\right)}, m \in \mathbb{N}_{0}.
\end{equation*}
Then the continuity of each {eigenvalue} $\lambda _{1 ,n}\left( p,q\right) $ in $p\in ({\mathcal{M}}_{0}(I,{\mathbb{R}}),w^{\ast })$ can be {proved by an argument similar to the one used} in the proof of (i).

Similarly, we can deduce that each {eigenvalue} $\lambda _{2 ,n}\left( p,q\right) $ is continuous in $p\in ({\mathcal{M}}_{0}(I,{\mathbb{R}}),w^{\ast })$. \qed

\begin{rem}
\label{lauc}
{Since the weak$^{\ast }$ topology is weaker than the strong topology induced by the norm $\Vert \cdot \Vert_{\mathbf{V}}$, it yields that for any fixed $q \in {\mathcal{M}}_{0}(I,{\mathbb{R}})$, $\xi =1,2$, the eigenvalue $\lambda_{\xi,n}$ is continuous in $p\in ({\mathcal{M}}_{0}(I,{\mathbb{R}}),\Vert \cdot \Vert_{\mathbf{V}})${, and for any fixed $p \in {\mathcal{M}}_{0}(I,{\mathbb{R}})$, $\xi =1,2$, the eigenvalue $\lambda_{\xi,n}$ is continuous in $q\in ({\mathcal{M}}_{0}(I,{\mathbb{R}}),\Vert \cdot \Vert_{\mathbf{V}})$.}}
\end{rem}

Finally, we deduce the differentiability of eigenvalues with respect to $p$, $q \in ({\mathcal{M}}_{0}(I,{\mathbb{R}}),\Vert \cdot \Vert _{\mathbf{V}})$.

\smallskip\smallskip\noindent\emph{Proof of Theorem \ref{frechet}.} (i) For any $p_{0}\in ({\mathcal{M}}_{0}(I,{\mathbb{R}}),\Vert \cdot \Vert_{\mathbf{V}})$, let
\begin{eqnarray*}
N_{1,p_{0}}\hh&:=&\hh\min \left\{ N\in {\mathbb{N}}\bigg |\left( 2N+1\right) \pi >\frac{9C_{\pi }}{4}e^{3\left(
3\Vert q \Vert_{\mathbf{V}}+\Vert p_{0} \Vert_{\mathbf{V}}  \right) },\right.\\
&&\left.2N\pi >\max\left\{ \frac{9C_{\pi }}{2}e^{3\left( 3\Vert q \Vert_{\mathbf{V}}+\Vert p_{0} \Vert_{\mathbf{V}} \right) },2\ln \frac{C_{\pi }}{2}, 2\sqrt{2}\ln \frac{C_{\pi }}{4}\right\} \right\}.
\end{eqnarray*}

\textit{Step 1.} For $\xi=1,2$, $\vert n\vert \geqslant N_{1,p_{0}}$, we first prove that the eigenfunction $E_{\xi ,n}\left( x,p,q\right)$ is {continuous in} $p \in ({\mathcal{M}}_{0}(I,{\mathbb{R}}),\Vert \cdot \Vert_{\mathbf{V}})$ {uniformly for $x\in I$}, {{i.e.}},
\begin{equation}
E_{\xi ,n}\left( \cdot,p,q\right) \rightarrow E_{\xi ,n}\left( \cdot,p_{0},q\right) \mbox{ as } \Vert p-p_{0}\Vert _{\mathbf{V}}\rightarrow 0  \label{Enp0}
\end{equation}
holds uniformly for $x \in I$. For $p$, $p_{0}\in ({\mathcal{M}}_{0}(I,{\mathbb{R}}),\Vert \cdot \Vert_{\mathbf{V}})$, denote
\begin{equation*}
M_{\xi }\left( \lambda_{\xi,n}(p,q) ,p,q\right) =\left(
\begin{array}{cc}
y_{1}(1,\lambda_{\xi,n}(p,q) ,p,q) & y_{2}(1,\lambda_{\xi,n}(p,q) ,p,q) \\
y_{1}^{\prime }(1,\lambda_{\xi,n}(p,q),p,q) & y_{2}^{\prime }(1,\lambda_{\xi,n}(p,q),p,q)+(-1)^{\xi}%
\end{array}%
\right)
\end{equation*}
and
\begin{equation}
M_{\xi }\left( \lambda_{\xi,n}(p_{0},q) ,p_{0},q\right) =\left(
\begin{array}{cc}
y_{1}(1,\lambda_{\xi,n}(p_{0},q) ,p_{0},q) & y_{2}(1,\lambda_{\xi,n}(p_{0},q) ,p_{0},q) \\
y_{1}^{\prime }(1,\lambda_{\xi,n}(p_{0},q),p_{0},q) & y_{2}^{\prime }(1,\lambda_{\xi,n}(p_{0},q),p_{0},q)+(-1)^{\xi}%
\end{array}%
\right) {.}\label{mp0}
\end{equation}
{T}hen according to {the fact}
\begin{equation*}
detM_{\xi }\left( \lambda_{\xi,n}(p,q) ,p,q\right)=detM_{\xi }\left( \lambda_{\xi,n}(p_{0},q) ,p_{0},q\right)=0,
\end{equation*}
we obtain
 \begin{equation*}
 0 \leqslant {\textrm{Rank}}M_{\xi }\left( \lambda_{\xi,n}(p,q) ,p,q\right)<2,\quad 0 \leqslant {\textrm{Rank}}M_{\xi }\left( \lambda_{\xi,n}(p_{0},q) ,p_{0},q\right)<2.
 \end{equation*}
Moreover, it follows from Corollary \ref{simple} that the $g$-multiplicity of each eigenvalue $\lambda _{\xi ,n}\left(p_{0},q\right)$, $\left\vert n\right\vert \geqslant N_{1,p_{0}}$, $\xi =1,2$, is one, then ${\textrm{Rank}}M_{\xi }\left( \lambda_{\xi,n}(p_{0},q) ,p_{0},q\right)=1${, {i.e.}, at least one entry of the matrix (\ref{mp0}) is nonzero.

\textbf{Case 1.} Suppose}
\begin{equation}
y_{1}(1,\lambda_{\xi,n}(p_{0},q) ,p_{0},q)\neq0, \left\vert n\right\vert \geqslant N_{1,p_{0}}, \xi =1,2, \label{neq0}
\end{equation}
and let
\begin{equation*}
a(p_{0},q):=y_{1}^{-1}(1,\lambda_{\xi,n}(p_{0},q) ,p_{0},q)y_{2}(1,\lambda_{\xi,n}(p_{0},q) ,p_{0},q),\quad b(p_{0},q):=-1{.}
\end{equation*}
{Then it is easy to see that}
\begin{equation*}
a(p_{0},q)y_{1}(1,\lambda_{\xi,n}(p_{0},q) ,p_{0},q)+b(p_{0},q)y_{2}(1,\lambda_{\xi,n}(p_{0},q) ,p_{0},q)=0.
\end{equation*}
Therefore, the eigenfunction corresponding to the eigenvalue $\lambda_{\xi,n}(p_{0},q)$, $\left\vert n\right\vert \geqslant N_{1,p_{0}}$, $\xi=1,2$, is
\begin{equation}
e(x,\lambda_{\xi,n}(p_{0},q),p_{0},q)\!\!=\!\!a(p_{0},q)y_{1}(x,\lambda_{\xi,n}(p_{0},q) ,p_{0},q)+b(p_{0},q)y_{2}(x,\lambda_{\xi,n}(p_{0},q) ,p_{0},q).\label{eigenfunction}
\end{equation}
According to {Remark \ref{lauc}},  {\PP}{\ref{uc2} $(ii)$} and {the fact} (\ref{neq0}), {it follows} that for each $n \geqslant N_{1,p_{0}}$, $\xi=1,2$, there exists a constant $\delta_{n} >0$ such that if $\Vert p-p_{0}\Vert_{\mathbf{V}} <\delta_{n}$, one has $y_{1}(1,\lambda_{\xi,n}(p,q) ,p,q)\neq0$. Let
\begin{equation*}
a(p,q):=y_{1}^{-1}(1,\lambda_{\xi,n}(p,q) ,p,q)y_{2}(1,\lambda_{\xi,n}(p,q) ,p,q),\quad b(p,q):=-1,
\end{equation*}
then
\begin{equation*}
a(p,q)y_{1}(1,\lambda_{\xi,n}(p,q) ,p,q)+b(p,q)y_{2}(1,\lambda_{\xi,n}(p,q) ,p,q)=0.
\end{equation*}
Therefore, when $\left\vert n\right\vert \geqslant N_{1,p_{0}}$, $\Vert p-p_{0}\Vert _{\mathbf{V}}<\delta_{n}$, the $g$-multiplicity of $\lambda_{\xi,n}(p,q)$, $\xi=1,2$, is one{,} and the corresponding eigenfunction is
\begin{equation*}
e(x,\lambda_{\xi,n}(p,q),p,q)=a(p,q)y_{1}(x,\lambda_{\xi,n}(p,q) ,p,q)+b(p,q)y_{2}(x,\lambda_{\xi,n}(p,q) ,p,q).
\end{equation*}
Additionally, from  {\PP}{\ref{uc2} $(ii)$}, one has (\ref{Enp0}) {holds in this case}.

{\textbf{Case 2.} Suppose
\begin{equation*}
y_{2}(1,\lambda_{\xi,n}(p_{0},q),p_{0},q)\neq 0,\quad \vert n\vert \geqslant N_{1,p_{0}},\quad \xi=1,2,
\end{equation*}
then we can define $a(p_{0},q)$, $b(p_{0},q)$ as follows:
\begin{equation*}
a(p_{0},q):=-1,\qquad b(p_{0},q):=y_{2}^{-1}(1,\lambda_{\xi,n}(p_{0},q),p_{0},q)y_{1}(1,\lambda_{\xi,n}(p_{0},q),p_{0},q).
\end{equation*}
It is easy to see that
\begin{equation*}
a(p_{0},q)y_{1}(1,\lambda_{\xi,n}(p_{0},q),p_{0},q)+b(p_{0},q)y_{2}(1,\lambda_{\xi,n}(p_{0},q),p_{0},q)=0.
\end{equation*}
Thus the eigenfunction corresponding to the eigenvalue $\lambda_{\xi,n}(p_{0},q)$ is (\ref{eigenfunction}), and then using the same argument as in the proof of \textbf{Case 1}, we can prove (\ref{Enp0}) in this case.}

{\textbf{Case 3.} Suppose
\begin{equation*}
y_{1}^{\prime}(1,\lambda_{\xi,n}(p_{0},q),p_{0},q)\neq 0,\quad \vert n\vert \geqslant N_{1,p_{0}},\quad \xi=1,2,
\end{equation*}
then we can define $a(p_{0},q)$, $b(p_{0},q)$ as follows:
\begin{equation*}
a(p_{0},q):=\left(y_{1}^{\prime}\right)^{-1}(1,\lambda_{\xi,n}(p_{0},q),p_{0},q)\left(y_{2}^{\prime}(1,\lambda_{\xi,n}(p_{0},q),p_{0},q)+(-1)^{\xi}\right),\qquad b(p_{0},q):=-1,
\end{equation*}
and
\begin{equation*}
a(p_{0},q)y_{1}^{\prime}(1,\lambda_{\xi,n}(p_{0},q),p_{0},q)+b(p_{0},q)\left(y_{2}^{\prime}(1,\lambda_{\xi,n}(p_{0},q),p_{0},q)+(-1)^{\xi}\right)=0.
\end{equation*}
Thus the eigenfunction corresponding to the eigenvalue $\lambda_{\xi,n}(p_{0},q)$ is (\ref{eigenfunction}), and then using the same argument as in the proof of \textbf{Case 1}, we can prove (\ref{Enp0}) in this case.}

{\textbf{Case 4.} Suppose
\begin{equation*}
y_{2}^{\prime}(1,\lambda_{\xi,n}(p_{0},q),p_{0},q)+(-1)^{\xi}\neq 0,\quad \vert n\vert \geqslant N_{1,p_{0}},\quad \xi=1,2,
\end{equation*}
then we can define $a(p_{0},q)$, $b(p_{0},q)$ as follows:
\begin{equation*}
a(p_{0},q):=-1,\qquad b(p_{0},q):=\left(y_{2}^{\prime}(1,\lambda_{\xi,n}(p_{0},q),p_{0},q)+(-1)^{\xi}\right)^{-1}y_{1}^{\prime}(1,\lambda_{\xi,n}(p_{0},q),p_{0},q),
\end{equation*}
and
\begin{equation*}
a(p_{0},q)y_{1}^{\prime}(1,\lambda_{\xi,n}(p_{0},q),p_{0},q)+b(p_{0},q)\left(y_{2}^{\prime}(1,\lambda_{\xi,n}(p_{0},q),p_{0},q)+(-1)^{\xi}\right)=0.
\end{equation*}
Thus the eigenfunction corresponding to the eigenvalue $\lambda_{\xi,n}(p_{0},q)$ is (\ref{eigenfunction}), and then using the same argument as in the proof of \textbf{Case 1}, we can prove (\ref{Enp0}) in this case.}

\textit{Step 2.} Now we deduce the Fr\'{e}chet derivatives of {the} eigenvalues $\lambda_{\xi,n}(p,q)$, $\xi =1,2$, $\left\vert n\right\vert \geqslant N_{1,p_{0}}$ at $p=p_{0}\in ({\mathcal{M}}_{0}(I,{\mathbb{R}}),\Vert\cdot \Vert _{\mathbf{V}})$. According to {the} equation (\ref{equation}), it follows that for $\nu_{p} \in ({\mathcal{M}}_{0}(I,{\mathbb{R}}),\Vert\cdot \Vert _{\mathbf{V}})$,
\begin{eqnarray*}
&&\hh(\lambda _{\xi ,n}\left( p_{0}+ \nu _{p},q\right)-\lambda _{\xi ,n}\left( p_{0},q\right))\int_{[0,1]}E_{\xi ,n}\left( x,p_{0}+\nu _{p},q\right)\bar{E}_{\xi ,n}\left( x,p_{0},q\right)\mathrm{d}x \\
&=&\hh\int_{[0,1]}\bar{E}_{\xi ,n}\left( x,p_{0},q\right)\left[E_{\xi ,n}\left( x,p_{0}+ \nu _{p},q\right)\left( i\mathrm{d}q\left( x\right)+\mathrm{d}p_{0}\left(x\right) +\mathrm{d}\nu _{p}\left( x\right)\right) \right.\\
& &\hh\left.+i\mathrm{d}\left( E_{\xi ,n}^{ \prime }\right) ^{\bullet }\left( x,p_{0}+ \nu_{p},q\right) +2iq\left( x\right) E_{\xi ,n}^{ \prime }\left( x,p_{0}+\nu_{p},q\right)\mathrm{d}x\right]-\left[\left( -i\mathrm{d}q\left( x\right)\right.\right.\\
& &\hh\left.\left.+\mathrm{d}p_{0}\left(x\right) \right)\bar{E}_{\xi ,n}\left( x,p_{0},q\right) -2iq\left( x\right) \bar{E}_{\xi ,n}^{ \prime }\left( x,p_{0},q\right)\mathrm{d}x\right.\\
& &\hh\left.-i\mathrm{d}\left( \bar{E}_{\xi ,n}^{ \prime }\right) ^{\bullet }\left( x,p_{0},q\right) \right]E_{\xi ,n}\left( x,p_{0}+\nu _{p},q\right).
\end{eqnarray*}
Using the integration by parts formula and the boundary conditions (BC)$_{\xi}$, one can deduce that
\begin{eqnarray*}
&&\hh(\lambda _{\xi ,n}\left( p_{0}+ \nu _{p},q\right)-\lambda _{\xi ,n}\left( p_{0},q\right))\int_{[0,1]}E_{\xi ,n}\left( x,p_{0}+\nu _{p},q\right)\bar{E}_{\xi ,n}\left( x,p_{0},q\right)\mathrm{d}x \\
&=&\hh\int_{[0,1]}\bar{E}_{\xi ,n}\left( x,p_{0},q\right)E_{\xi ,n}\left( x,p_{0}+ \nu _{p},q\right)\mathrm{d}\nu _{p}.
\end{eqnarray*}
Dividing both sides by $\nu_{p}$, letting $\Vert \nu_{p} \Vert_{\mathbf{V}}\rightarrow 0${,} and using the {statement in} (\ref{Enp0}), one has
\begin{equation*}
\partial _{p}\lambda _{\xi ,n}\left( p_{0},q\right)=\left\vert E_{\xi ,n}\left( x,p_{0},q\right) \right\vert ^{2}.
\end{equation*}
Since $\left\vert E_{\xi ,n}\left( x,p_{0},q\right) \right\vert ^{2}\in {\mathcal{C}%
}(I,{\mathbb{R}})$, it follows that each $\partial _{p}\lambda _{\xi ,n}\left(
p_{0},q\right)$, $\left\vert n\right\vert \geqslant N_{1,p_{0}}$, $\xi =1,2$, is a bounded linear functional of $({\mathcal{M}}_{0}(I,{\mathbb{%
R}}),\Vert \cdot \Vert _{\mathbf{V}})$, that is, $\partial _{p}\lambda
_{\xi ,n}\left( p_{0},q\right) \in ({\mathcal{M}}_{0}(I,{\mathbb{R}}),\Vert
\cdot \Vert _{\mathbf{V}})^{\ast }\cong ({\mathcal{C}}(I,{\mathbb{R}}),\Vert
\cdot \Vert _{\infty })^{\ast \ast }$.

(ii) For $ q_{0}\in ({\mathcal{M}}_{0}(I,{\mathbb{R}}),\Vert \cdot \Vert_{\mathbf{V}})$, let
\begin{eqnarray*}
N_{2,q_{0}}\hh&:=&\hh\min \left\{ N\in {\mathbb{N}}\bigg |\left( 2N+1\right) \pi >\frac{9C_{\pi }}{4}e^{3\left(3 \Vert q_{0} \Vert_{\mathbf{V}}+\Vert p \Vert_{\mathbf{V}} \right) },\right.\\
&&\left.2N\pi >\max
\left\{ \frac{9C_{\pi }}{2}e^{3\left( 3 \Vert q_{0} \Vert_{\mathbf{V}}+\Vert p \Vert_{\mathbf{V}} \right) },2\ln \frac{C_{\pi }}{2}, 2\sqrt{2}\ln \frac{C_{\pi }}{4}\right\} \right\}.
\end{eqnarray*}

\textit{Step 1.} For $\xi=1,2$, $\vert n\vert \geqslant N_{2,q_{0}}$, using the same {argument} as {in} the proof of Theorem \ref{frechet} (i), we can prove that the eigenfunction $E_{\xi ,n}\left( x,p,q\right)$ is {continuous} in $q \in ({\mathcal{M}}_{0}(I,{\mathbb{R}}),\Vert \cdot \Vert_{\mathbf{V}})$ {for} $x\in I$, {{i.e.}},
\begin{equation}
E_{\xi ,n}\left( \cdot,p,q\right) \rightarrow E_{\xi ,n}\left( \cdot,p,q_{0}\right) \mbox{ as } \Vert q-q_{0}\Vert _{\mathbf{V}}\rightarrow 0  \label{Enq0}
\end{equation}
holds uniformly for $x \in I$.

\textit{Step 2.} Now we deduce the Fr\'{e}chet derivatives of {the} eigenvalues $\lambda_{\xi,n}(p,q)$, $\xi =1,2$, $\left\vert n\right\vert \geqslant N_{2,q_{0}}$ at $q=q_{0}\in ({\mathcal{M}}_{0}(I,{\mathbb{R}}),\Vert\cdot \Vert _{\mathbf{V}})$.
{According} to {the} equation (\ref{equation}), it follows that for $\nu_{q} \in ({\mathcal{M}}_{0}(I,{\mathbb{R}}),\Vert\cdot \Vert _{\mathbf{V}})$,
\begin{eqnarray*}
&&\hh(\lambda _{\xi ,n}\left( p,q_{0}+ \nu _{q}\right)-\lambda _{\xi ,n}\left( p,q_{0}\right))\int_{[0,1]}E_{\xi ,n}\left( x,p,q_{0}+\nu _{q}\right)\bar{E}_{\xi ,n}\left( x,p,q_{0}\right)\mathrm{d}x \\
&=&\hh\int_{[0,1]}\bar{E}_{\xi ,n}\left( x,p,q_{0}\right)\Bigl[E_{\xi ,n}\left( x,p,q_{0}+ \nu _{{q}}\right)\left( i\mathrm{d}q_{0}\left( x\right)+i\mathrm{d}\nu _{q}\left( x\right)+\mathrm{d}p\left(x\right) \right) \Bigr.\\
& &\hh\left.+i\mathrm{d}\left( E_{\xi ,n}^{ \prime }\right) ^{\bullet }\left( x,p,q_{0}+ \nu_{q}\right) +2i(q_{0}\left( x\right)+\nu_{q}(x)) E_{\xi ,n}^{ \prime }\left( x,p,q_{0}+\nu_{q}\right)\mathrm{d}x\right]\\
& &\hh-E_{\xi ,n}\left( x,p,q_{0}+\nu _{q}\right)\left[ -i\mathrm{d}\left( \bar{E}_{\xi ,n}^{ \prime }\right) ^{\bullet }\left( x,p,q_{0}\right) -2iq_{0}\left( x\right) \bar{E}_{\xi ,n}^{ \prime }\left( x,p,q_{0}\right)\mathrm{d}x\right.\\
& &\hh\Bigl.+\bar{E}_{\xi ,n}\left( x,p,q_{0}\right)\left( -i\mathrm{d}q_{0}\left( x\right)+\mathrm{d}p\left(x\right) \right)\Bigr].
\end{eqnarray*}
{Using the integration by parts formula, we obtain
\begin{eqnarray*}
&&\!\!\!\!(\lambda _{\xi ,n}\left( p,q_{0}+ \nu _{q}\right)-\lambda _{\xi ,n}\left( p,q_{0}\right))\int_{[0,1]}E_{\xi ,n}\left( x,p,q_{0}+\nu _{q}\right)\bar{E}_{\xi ,n}\left( x,p,q_{0}\right)\mathrm{d}x \\
&=&\!\!\!\! i\left[(E^{\prime})^{\bullet}_{\xi ,n}\left( x,p,q_{0}+\nu _{q}\right)\bar{E}_{\xi ,n}\left( x,p,q_{0}\right)-E^{\prime}_{\xi ,n}\left( x,p,q_{0}+\nu _{q}\right)\bar{E}^{\prime}_{\xi ,n}\left( x,p,q_{0}\right)\right.\\
&&\!\!\!\! \left.+E_{\xi ,n}\left( x,p,q_{0}+\nu _{q}\right)(\bar{E}_{\xi ,n}^{\prime})^{\bullet}\left( x,p,q_{0}\right)+2q_{0}(x)E_{\xi ,n}\left( x,p,q_{0}+\nu _{q}\right)\bar{E}_{\xi ,n}\left( x,p,q_{0}\right)\right.\\
&&\!\!\!\! \left.+\nu _{q}(x)E_{\xi ,n}\left( x,p,q_{0}+\nu _{q}\right)\bar{E}_{\xi ,n}\left( x,p,q_{0}\right)\right]\vert_{x=0}^{x=1+}\\
&&\!\!\!\!+\int_{[0,1]}i[\bar{E}_{\xi ,n}\left( x,p,q_{0}\right)E^{\prime}_{\xi ,n}\left( x,p,q_{0}+\nu _{q}\right)-\bar{E}^{\prime}_{\xi ,n}\left( x,p,q_{0}\right)E_{\xi ,n}\left( x,p,q_{0}+\nu _{q}\right)]\nu _{q}(x)\mathrm{d}x.
\end{eqnarray*}
Then according to {the fact} $q_{0}(0)=\nu _{q}(0)=0$ and the boundary conditions (BC)$_{\xi}$, we obtain
\begin{eqnarray*}
&&\!\!\!\!(\lambda _{\xi ,n}\left( p,q_{0}+ \nu _{q}\right)-\lambda _{\xi ,n}\left( p,q_{0}\right))\int_{[0,1]}E_{\xi ,n}\left( x,p,q_{0}+\nu _{q}\right)\bar{E}_{\xi ,n}\left( x,p,q_{0}\right)\mathrm{d}x \\
&=&\!\!\!\! \int_{[0,1]}i[\bar{E}_{\xi ,n}\left( x,p,q_{0}\right)E^{\prime}_{\xi ,n}\left( x,p,q_{0}+\nu _{q}\right)-\bar{E}^{\prime}_{\xi ,n}\left( x,p,q_{0}\right)E_{\xi ,n}\left( x,p,q_{0}+\nu _{q}\right)]\nu _{q}(x)\mathrm{d}x.
\end{eqnarray*}
Using the fact {(\ref{Enq0})} we get
\begin{equation*}
(\lambda _{\xi ,n}\left( p,q_{0}+ \nu _{q}\right)\!-\!\lambda _{\xi ,n}\left( p,q_{0}\right))(1+o(1))=i\!\!\int_{[0,1]}\!\!\left[ E_{\xi ,n}\left( x,p,q_{0}\right) ,\bar{E}_{\xi ,n}\left( x,p,q_{0}\right) \right] \nu _{q}(x)\mathrm{d}x+o(\Vert \nu _{q}\Vert_{\mathbf{V}})
\end{equation*}
as $\nu _{q}\rightarrow 0$ in $({\mathcal{M}}_{0}(I,{\mathbb{R}}),\Vert \cdot \Vert_{\mathbf{V}})$. Consequently,
\begin{eqnarray*}
\lambda _{\xi ,n}\left( p,q_{0}+ \nu _{q}\right)-\lambda _{\xi ,n}\left( p,q_{0}\right)\!\!\!\!&=&\!\!\!\!\Bigg(i\int_{[0,1]}\left[ E_{\xi ,n}\left( x,p,q_{0}\right) ,\bar{E}_{\xi ,n}\left( x,p,q_{0}\right) \right] \nu _{q}(x)\mathrm{d}x\\
&&\!\!\!\!+o(\Vert \nu _{q}\Vert_{\mathbf{V}})\Bigg)(1+o(1))^{-1}\\
\!\!\!\!&=&\!\!\!\!i\int_{[0,1]}\left[ E_{\xi ,n}\left( x,p,q_{0}\right) ,\bar{E}_{\xi ,n}\left( x,p,q_{0}\right) \right] \nu _{q}(x)\mathrm{d}x+o(\Vert \nu _{q}\Vert_{\mathbf{V}})
\end{eqnarray*}
as $\nu _{q}\rightarrow 0$ in $({\mathcal{M}}_{0}(I,{\mathbb{R}}),\Vert \cdot \Vert_{\mathbf{V}})$. Hence,}
\begin{equation*}
\partial _{q}\lambda _{\xi ,n}\left( p,q_{0}\right)=i\left[ E_{\xi ,n}\left( x,p,q_{0}\right) ,\bar{E}_{\xi ,n}\left( x,p,q{_{0}}\right) \right] .
\end{equation*}
{It is easy to see that} $i\left[ E_{\xi ,n}\left( x,p,q_{0}\right) ,\bar{E}_{\xi ,n}\left( x,p,q{_{0}}\right) \right]\in {\mathcal{C}}(I,{\mathbb{R}})${. H}ence each
$\partial _{q}\lambda _{\xi ,n}\left( p,q_{0}\right) $, $\xi =1,2$, $\left\vert n\right\vert \geqslant N_{2,q_{0}}$, is a bounded linear
functional of $({\mathcal{M}}_{0}(I,{\mathbb{R}}),\Vert \cdot \Vert _{%
\mathbf{V}})$, that is, $\partial _{q}\lambda _{\xi ,n}\left( p,q_{0}\right)
\in ({\mathcal{M}}_{0}(I,{\mathbb{R}}),\Vert \cdot \Vert _{\mathbf{V}%
})^{\ast }\cong ({\mathcal{C}}(I,{\mathbb{R}}),\Vert \cdot \Vert _{\infty
})^{\ast \ast }$. \qed

\section*{Acknowledgements}
The research is supported by the National Natural Science Foundation of China (Grant No. 11601372); the Science and Technology Research Project of Higher Education in Hebei Province (Grant No. QN2017044).

\end{document}